\documentclass[11pt]{amsart}

\usepackage[utf8]{inputenc}
\usepackage{fullpage}
\usepackage[colorlinks,linkcolor=blue,citecolor=blue]{hyperref}
\usepackage{graphicx}
\usepackage{amsmath}
\usepackage{amsthm}
\usepackage{amsfonts}
\usepackage{amssymb}
\usepackage{bbm}
\usepackage{mathtools}
\usepackage{extarrows}
\usepackage{array}   % for \newcolumntype macro
\usepackage{booktabs, longtable}
\newcolumntype{L}{>{$}l<{$}}
\usepackage{tikz-cd}
\usetikzlibrary{arrows.meta}
\usepackage{float}
\usepackage{scalerel}

\usepackage{amsmath,tikz}
\usetikzlibrary{
  knots,
  hobby,
  decorations.pathreplacing,
  shapes.geometric,
  calc
}
\usetikzlibrary{arrows.meta,decorations.markings}

\usepackage{amsmath}
\usetikzlibrary{arrows.meta,decorations.markings,calc}
\usepackage{mathrsfs}
\usetikzlibrary{arrows}

\newtheorem{Thm}{Theorem}[section]

\newtheorem{Cor}[Thm]{Corollary}
\newtheorem{prop}[Thm]{Proposition}
\newtheorem{lem}[Thm]{Lemma}

\newtheorem{defn}[Thm]{Definition}

\theoremstyle{definition}

\theoremstyle{remark}

\newtheorem{rmk}{Remark}[section]
\newcommand{\one}{\mathbf{1}}

\usepackage{eucal}
\usepackage{graphicx}
\usepackage{tabularx}
\usepackage{multirow}
\usepackage[all,knot]{xy}

\DeclareMathOperator{\Rep}{Rep}

\DeclareMathOperator{\rev}{rev}

\DeclareMathOperator{\PSU}{PSU}

\DeclareMathOperator{\End}{End}

\DeclareMathOperator{\Tr}{Tr}

\DeclareMathOperator{\Irr}{Irr}
\DeclareMathOperator{\id}{id}
\DeclareMathOperator{\Id}{Id}

\DeclareMathOperator{\Obj}{Obj}
\newcommand{\RNum}[1]{\uppercase\expandafter{\romannumeral #1\relax}}

\newcommand{\BZ}{\mathbb{Z}}
\newcommand{\Zn}[1]{\BZ/{#1}\BZ}
\DeclareMathOperator{\Aut}{Aut}

\usepackage{graphicx}
\graphicspath{ {./graph/} }

\newcommand{\BB}{\mathcal{B}}
\newcommand{\BC}{\mathbb{C}}

\renewcommand{\AA}{\mathcal{A}}
\newcommand{\CC}{{\mathcal{C}}}
\newcommand{\DD}{\mathcal{D}}

\newcommand{\ZZ}{\mathcal{Z}}
\newcommand{\EE}{\mathcal{E}}

\newcommand{\Hom}{\operatorname{Hom}}

\newcommand{\FPdim}{\operatorname{FPdim}}

\newcommand{\tr}{\mathrm{tr}}

\newcommand{\Inv}{\operatorname{Inv}}

\newcommand{\OverlapCoupons}[2]{%
  \mathord{\vcenter{\hbox{%
    \begin{tikzpicture}[
      x=5.5mm,y=5.5mm,draw=blue,
      line width=.7pt,line cap=round]
      \foreach \x in {-2.6,-.65,0,.65,2.6}{
        \draw (\x,0)--(\x,4.2);
      }
      \foreach \y in {3.75}{
        \node at (-1.625,\y) {$\cdots$};
        \node at ( 1.625,\y) {$\cdots$};
      }
      \draw[fill=white]
        (-3,2.5) rectangle (.35,3.2);
      \node at (-1.325,2.85) {$#1$};

      \draw[fill=white]
        (-.35,1.0) rectangle (3,1.7);
      \node at (1.325,1.35) {$#2$};

    \end{tikzpicture}%
  }}}%
}

\tikzset{
  couponcalc/.style={
    x=5mm,y=5mm,draw=blue,line width=.65pt,
    line cap=butt,line join=miter
  }
}

\newcommand{\CouponCableLayout}[1]{%
  \ifnum#1=0
    \def\CouponCalcStrands{-2.8,-2.1,0,.7,2.8}%
  \else
    \def\CouponCalcStrands{-2.8,-.7,0,.7,2.8}%
  \fi
}
\newcommand{\CouponCableDots}[2]{%
  \ifnum#1=0 \def\CouponDotsLeft{-1.05}%
  \else \def\CouponDotsLeft{-1.75}\fi
  \node[font=\scriptsize,inner sep=0pt]
    at (\CouponDotsLeft,#2) {$\cdots$};
  \node[font=\scriptsize,inner sep=0pt]
    at (1.75,#2) {$\cdots$};
}
\newcommand{\DrawCouponOp}[2]{%
  \ifnum#1=0
    \def\CouponCalcLayout{0}%
  \else
    \ifdefined\CouponCalcLayout\else\def\CouponCalcLayout{1}\fi
  \fi
  \CouponCableLayout{\CouponCalcLayout}
  \foreach \x in \CouponCalcStrands{
    \draw[#2,arrows=-] (\x,0)--(\x,1.8);
  }
  \ifnum#1=0
    \draw[fill=white] (-3.2,.4) rectangle (.35,1.4);
    \node[font=\footnotesize] at (-1.425,.9) {$f$};
  \else
    \draw[fill=white] (-.35,.4) rectangle (3.2,1.4);
    \node[font=\footnotesize] at (1.425,.9) {$f$};
  \fi
  \ifdefined\CouponCalcOmitDots\else
    \CouponCableDots{\CouponCalcLayout}{1.55}
  \fi
}

\newcommand{\DrawBlockSwap}[1]{%
  \draw[#1,arrows=-]
    (0,0) .. controls (0,.8) and (-2.8,2) .. (-2.8,2.8);
  \draw[#1,arrows=-]
    (.7,0) .. controls (.7,.8) and (-2.1,2) .. (-2.1,2.8);
  \draw[#1,arrows=-]
    (2.8,0) .. controls (2.8,.8) and (0,2) .. (0,2.8);
  \draw[#1,arrows=-,preaction={
    draw=white,line width=2.6pt,line cap=butt}]
    (-2.8,0) .. controls (-2.8,.8) and (.7,2) .. (.7,2.8);
  \draw[#1,arrows=-,preaction={
    draw=white,line width=2.6pt,line cap=butt}]
    (-.7,0) .. controls (-.7,.8) and (2.8,2) .. (2.8,2.8);
  \ifdefined\CouponCalcOmitDots\else
    \CouponCableDots{0}{2.55}
  \fi
}
\newcommand{\CouponOp}[1]{%
  \mathord{\vcenter{\hbox{%
    \begin{tikzpicture}[couponcalc]
      \DrawCouponOp{#1}{}
    \end{tikzpicture}%
  }}}%
}

\newcommand{\CouponSlide}[1]{%
  \mathord{\vcenter{\hbox{%
    \begin{tikzpicture}[couponcalc]
      \def\CouponCalcOmitDots{}%
      \ifnum#1=2
        \DrawCouponOp{1}{}
        \begin{scope}[shift={(0,1.8)}]
          \DrawBlockSwap{}
        \end{scope}
      \else
        \DrawBlockSwap{}
        \begin{scope}[shift={(0,2.8)}]
          \def\CouponCalcLayout{0}%
          \ifnum#1=1
            \DrawCouponOp{0}{}
          \else
            \DrawCouponOp{1}{}
          \fi
        \end{scope}
      \fi
      \CouponCableDots{0}{4.35}
    \end{tikzpicture}%
  }}}%
}

\tikzset{
  tracecalc/.style={
    x=5mm,y=5mm,draw=blue,line width=.65pt,
    line cap=round,line join=round
  },
  tracecoupon/.style={
    draw=blue,fill=white,rectangle,rounded corners=0pt,
    line join=miter,minimum width=6.4mm,minimum height=3.6mm,
    inner xsep=.7pt,inner ysep=.5pt,font=\footnotesize
  }
}
\newcommand{\TraceReturn}[4]{%
  \draw (#1,#3)
    .. controls (#1,{#3+#4}) and (#2,{#3+#4}) .. (#2,#3)
    --(#2,0)
    .. controls (#2,-#4) and (#1,-#4) .. (#1,0);
  \draw (#2,{#3/2+.32})--(#2,{#3/2-.32});
}
\newcommand{\TraceCableClosure}[1]{%
  \TraceReturn{0}{4.9}{#1}{1.35}
  \TraceReturn{2}{3.3}{#1}{.6}
  \node[font=\scriptsize,inner sep=0pt] at (4.1,{#1/2}) {$\cdots$};
  \node[anchor=west,inner sep=2pt,font=\small]
    at (4.9,{#1/2}) {$X$};
}
\newcommand{\TraceXStrands}[2]{%
  \foreach \x in {0,2}{\draw (\x,#1)--(\x,#2);}
}
\newcommand{\TraceDots}[1]{%
  \node[font=\scriptsize,inner sep=0pt] at (1,#1) {$\cdots$};
}
\newcommand{\TraceCoupon}[2]{%
  \node[tracecoupon,minimum width=13.4mm] at (1,#1) {$#2$};
}

\newcommand{\TraceDiagram}[1]{%
  \mathord{\vcenter{\hbox{%
    \begin{tikzpicture}[tracecalc]
      \ifcase#1
        \TraceCableClosure{1.8}
        \TraceXStrands{0}{1.8}
        \TraceDots{.25}\TraceDots{1.55}
        \TraceCoupon{.9}{T}
      \or
        \TraceCableClosure{4.2}
        \TraceXStrands{0}{.6}
        \TraceXStrands{2.1}{4.2}
        \draw (1,.6)--(1,2.1);
        \TraceDots{2.85}
        \node[anchor=west,inner sep=2pt,font=\scriptsize]
          at (1,1.35) {$Y$};
        \TraceCoupon{.6}{p_a}
        \TraceCoupon{2.1}{i_a}
        \TraceCoupon{3.6}{T}
      \or
        \TraceReturn{1}{3.7}{4.2}{.9}
        \draw (1,0)--(1,.6) (1,3.6)--(1,4.2);
        \TraceXStrands{.6}{3.6}
        \TraceDots{1.35}\TraceDots{2.85}
        \node[anchor=east,inner sep=2pt,font=\scriptsize]
          at (0,1.35) {$X$};
        \node[anchor=west,inner sep=2pt,font=\small]
          at (3.7,2.1) {$Y$};
        \TraceCoupon{.6}{i_a}
        \TraceCoupon{2.1}{T}
        \TraceCoupon{3.6}{p_a}
      \or
        \TraceReturn{0}{2}{1.2}{.8}
        \draw (0,0)--(0,1.2);
        \node[anchor=west,inner sep=2pt,font=\small]
          at (2,.6) {$Y$};
      \fi
    \end{tikzpicture}%
  }}}%
}

\tikzset{
  crossingshift/.style={
    x=5.5mm,y=4.2mm,draw=blue,line width=.7pt,
    line cap=round,line join=round
  }
}

\newcommand{\ShiftCoordinateInLayout}[3]{%
  \ifdefined\GeneralShiftCabling
    \pgfmathsetmacro{#3}{{-3.85,-1.75,-1.05,1.05,1.75,3.85}[#2]}%
  \else
    \ifnum#1=0
      \pgfmathsetmacro{#3}{{-3.15,-1.05,-.35,.35,2.45,3.15}[#2]}%
    \else
      \pgfmathsetmacro{#3}{{-3.15,-2.45,-.35,.35,1.05,3.15}[#2]}%
    \fi
  \fi
}
\newcommand{\ShiftCoordinate}[2]{%
  \ShiftCoordinateInLayout{\ShiftDrawingLayout}{#1}{#2}%
}
\newcommand{\DrawShiftVertical}[2]{%
  \foreach \j in {0,...,5}{
    \ShiftCoordinate{\j}{\ShiftX}%
    \draw (\ShiftX,{#1})--(\ShiftX,{#2});
  }
}
\newcommand{\DrawShiftDots}[2]{%
  \ifnum#1=0
    \def\ShiftDotList{-2.1,1.4}%
  \else
    \def\ShiftDotList{-1.4,2.1}%
  \fi
  \foreach \x in \ShiftDotList{
    \node[font=\scriptsize,inner sep=0pt] at (\x,{#2}) {$\cdots$};
  }
}
\newcommand{\DrawShiftPath}[5]{%
  \draw[#5,arrows=-]
    (#1,{#3})
    .. controls (#1,{(#3)+.30*(#4)}) and (#2,{(#3)+.70*(#4)})
    .. (#2,{(#3)+(#4)});
}
\newcommand{\DrawShiftPair}[4]{%
  \ShiftCoordinate{#1}{\ShiftLeft}%
  \pgfmathtruncatemacro{\ShiftNext}{#1+1}%
  \ShiftCoordinate{\ShiftNext}{\ShiftRight}%
  \foreach \j in {0,...,5}{
    \ifnum\j=#1\relax\else
      \ifnum\j=\ShiftNext\relax\else
        \ShiftCoordinate{\j}{\ShiftX}%
        \draw (\ShiftX,#3)--(\ShiftX,{#3+#4});
      \fi
    \fi
  }
  \ifnum#2=1
    \DrawShiftPath{\ShiftRight}{\ShiftLeft}{#3}{#4}{}
    \DrawShiftPath{\ShiftLeft}{\ShiftRight}{#3}{#4}{
      preaction={draw=white,line width=2.8pt,line cap=butt}}
  \else
    \DrawShiftPath{\ShiftLeft}{\ShiftRight}{#3}{#4}{}
    \DrawShiftPath{\ShiftRight}{\ShiftLeft}{#3}{#4}{
      preaction={draw=white,line width=2.8pt,line cap=butt}}
  \fi
}
\newcommand{\DrawShiftBlock}[2]{%
  \ifnum#1=1
    \foreach \j in {4,5}{
      \pgfmathtruncatemacro{\ShiftDest}{\j-4}%
      \ShiftCoordinateInLayout{1}{\j}{\ShiftFrom}%
      \ShiftCoordinateInLayout{0}{\ShiftDest}{\ShiftTo}%
      \DrawShiftPath{\ShiftFrom}{\ShiftTo}{#2}{3}{}
    }
    \foreach \j in {0,...,3}{
      \pgfmathtruncatemacro{\ShiftDest}{\j+2}%
      \ShiftCoordinateInLayout{1}{\j}{\ShiftFrom}%
      \ShiftCoordinateInLayout{0}{\ShiftDest}{\ShiftTo}%
      \DrawShiftPath{\ShiftFrom}{\ShiftTo}{#2}{3}{
        preaction={draw=white,line width=2.8pt,line cap=butt}}
    }
    \def\ShiftDrawingLayout{0}%
  \else
    \foreach \j in {0,1}{
      \pgfmathtruncatemacro{\ShiftDest}{\j+4}%
      \ShiftCoordinateInLayout{0}{\j}{\ShiftFrom}%
      \ShiftCoordinateInLayout{1}{\ShiftDest}{\ShiftTo}%
      \DrawShiftPath{\ShiftFrom}{\ShiftTo}{#2}{3}{}
    }
    \foreach \j in {2,...,5}{
      \pgfmathtruncatemacro{\ShiftDest}{\j-2}%
      \ShiftCoordinateInLayout{0}{\j}{\ShiftFrom}%
      \ShiftCoordinateInLayout{1}{\ShiftDest}{\ShiftTo}%
      \DrawShiftPath{\ShiftFrom}{\ShiftTo}{#2}{3}{
        preaction={draw=white,line width=2.8pt,line cap=butt}}
    }
    \def\ShiftDrawingLayout{1}%
  \fi
}
\newcommand{\DrawShiftTwist}[2]{%
  \def\ShiftDrawingLayout{0}%
  \DrawShiftVertical{#2}{#2+.15}
  \DrawShiftPair{0}{#1}{#2+.15}{1}
  \DrawShiftVertical{#2+1.15}{#2+1.45}
  \DrawShiftPair{0}{#1}{#2+1.45}{1}
  \DrawShiftVertical{#2+2.45}{#2+2.6}
}
\newcommand{\DrawShiftH}[2]{%
  \ifnum#1=1
    \DrawShiftBlock{1}{#2}
    \DrawShiftTwist{1}{#2+3}
  \else
    \DrawShiftTwist{-1}{#2}
    \DrawShiftBlock{-1}{#2+2.6}
  \fi
}
\newcommand{\DrawShiftEndpoints}[2]{%
  \def\ShiftDrawingLayout{#1}%
  \foreach \j in {0,...,5}{
    \ShiftCoordinate{\j}{\ShiftX}%
    \draw (\ShiftX,{#2-.35})--(\ShiftX,#2);
  }
  \DrawShiftDots{#1}{#2-.18}
}
\newcommand{\DrawShiftCrossing}[3]{%
  \ifnum#1=1
    \DrawShiftPair{0}{1}{#2}{#3-#2}
  \else\ifnum#1=4
    \DrawShiftPair{4}{1}{#2}{#3-#2}
  \else
    \DrawShiftPair{2}{1}{#2}{#3-#2}
  \fi\fi
}
\newcommand{\CrossingShiftOperator}[2]{%
  \mathord{\vcenter{\hbox{%
    \begin{tikzpicture}[crossingshift]
      \ifnum#2=0
        \ifnum#1=2 \def\ShiftLayout{1}\else\def\ShiftLayout{0}\fi
        \def\ShiftDrawingLayout{\ShiftLayout}%
        \DrawShiftVertical{0}{1}
        \DrawShiftCrossing{#1}{1}{2.6}
        \DrawShiftVertical{2.6}{3.6}
        \DrawShiftEndpoints{\ShiftLayout}{3.6}
      \else
        \ifnum#2=1
          \def\ShiftLayout{0}%
          \def\ShiftDrawingLayout{0}%
          \DrawShiftVertical{0}{.4}
          \DrawShiftH{-1}{.4}
          \DrawShiftPair{0}{1}{6}{1.4}
          \DrawShiftH{1}{7.4}
        \else
          \def\ShiftLayout{1}%
          \def\ShiftDrawingLayout{1}%
          \DrawShiftVertical{0}{.4}
          \DrawShiftH{1}{.4}
          \DrawShiftPair{4}{1}{6}{1.4}
          \DrawShiftH{-1}{7.4}
        \fi
        \DrawShiftVertical{13}{13.4}
        \DrawShiftEndpoints{\ShiftLayout}{13.4}
      \fi
    \end{tikzpicture}%
  }}}%
}

\newcommand{\DrawGeneralShiftDots}[1]{%
  \foreach \x in {-2.8,0,2.8}{
    \node[font=\scriptsize,inner sep=0pt]
      at (\x,{#1}) {$\cdots$};
  }
}
\newcommand{\DrawGeneralShiftEndpoints}[1]{%
  \foreach \j in {0,...,5}{
    \ShiftCoordinate{\j}{\ShiftX}%
    \draw (\ShiftX,{#1-.35})--(\ShiftX,#1);
  }
  \DrawGeneralShiftDots{#1-.18}
}
\newcommand{\GeneralCrossingShiftOperator}[1]{%
  \mathord{\vcenter{\hbox{%
    \begin{tikzpicture}[crossingshift,x=4.5mm]
      \def\GeneralShiftCabling{}%
      \def\ShiftDrawingLayout{0}%
      \DrawGeneralShiftDots{.25}
      \ifnum#1=1
        \DrawShiftVertical{0}{.6}
        \DrawShiftBlock{-1}{.6}
        \DrawShiftPair{1}{1}{3.6}{1.6}
        \DrawGeneralShiftDots{4.4}
        \DrawShiftBlock{1}{5.2}
        \DrawShiftVertical{8.2}{8.8}
        \DrawGeneralShiftEndpoints{8.8}
      \else
        \DrawShiftVertical{0}{1}
        \DrawShiftPair{3}{1}{1}{1.6}
        \DrawShiftVertical{2.6}{3.6}
        \DrawGeneralShiftEndpoints{3.6}
      \fi
    \end{tikzpicture}%
  }}}%
}

\newcommand{\DrawTransportCoupon}[2]{%
  \DrawShiftVertical{#2}{(#2)+1.8}
  \ifnum#1=0
    \ShiftCoordinate{0}{\TransportFirst}%
    \ShiftCoordinate{3}{\TransportLast}%
  \else
    \ShiftCoordinate{2}{\TransportFirst}%
    \ShiftCoordinate{5}{\TransportLast}%
  \fi
  \draw[fill=white]
    ({\TransportFirst-.28},{(#2)+.225})
    rectangle ({\TransportLast+.28},{(#2)+1.575});
  \node at ({(\TransportFirst+\TransportLast)/2},{(#2)+.9}) {$f$};
}
\newcommand{\DrawTransportEnds}[2]{%
  \def\ShiftDrawingLayout{#1}%
  \foreach \j in {0,...,5}{
    \ShiftCoordinate{\j}{\ShiftX}%
    \draw (\ShiftX,{(#2)-.35})--(\ShiftX,#2);
  }
  \DrawShiftDots{#1}{(#2)-.18}
}
\newcommand{\TransportCommutationDiagram}[2]{%
  \mathord{\vcenter{\hbox{%
    \begin{tikzpicture}[crossingshift,x=4.05mm,y=3.6mm]
      \ifnum#1=1
        \def\TransportOutside{0}%
        \def\TransportOuterCoupon{1}%
        \def\TransportInnerCoupon{0}%
        \def\TransportInnerPair{0}%
        \def\TransportBottomSign{-1}%
        \def\TransportTopSign{1}%
        \def\TransportInside{1}%
      \else
        \def\TransportOutside{1}%
        \def\TransportOuterCoupon{0}%
        \def\TransportInnerCoupon{1}%
        \def\TransportInnerPair{4}%
        \def\TransportBottomSign{1}%
        \def\TransportTopSign{-1}%
        \def\TransportInside{0}%
      \fi
      \def\ShiftDrawingLayout{\TransportOutside}%
      \DrawShiftDots{\TransportOutside}{.2}
      \DrawShiftVertical{0}{.4}
      \ifnum#2=0
        \DrawShiftPair{2}{1}{.4}{1.6}
        \DrawTransportCoupon{\TransportOuterCoupon}{2}
        \DrawShiftVertical{3.8}{4.2}
        \DrawTransportEnds{\TransportOutside}{4.2}
      \else\ifnum#2=3
        \DrawTransportCoupon{\TransportOuterCoupon}{.4}
        \DrawShiftPair{2}{1}{2.2}{1.6}
        \DrawShiftVertical{3.8}{4.2}
        \DrawTransportEnds{\TransportOutside}{4.2}
      \else
        \DrawShiftBlock{\TransportBottomSign}{.4}
        \ifnum#2=1
          \DrawShiftPair{\TransportInnerPair}{1}{3.4}{1.6}
          \DrawTransportCoupon{\TransportInnerCoupon}{5}
        \else
          \DrawTransportCoupon{\TransportInnerCoupon}{3.4}
          \DrawShiftPair{\TransportInnerPair}{1}{5.2}{1.6}
        \fi
        \DrawShiftBlock{\TransportTopSign}{6.8}
        \DrawShiftVertical{9.8}{10.2}
        \DrawTransportEnds{\TransportOutside}{10.2}
      \fi\fi
    \end{tikzpicture}%
  }}}%
}
\newcommand{\CommuteDiagram}[1]{\TransportCommutationDiagram{-1}{#1}}
\newcommand{\ForwardCommuteDiagram}[1]{\TransportCommutationDiagram{1}{#1}}

\newcommand{\LabelledCouponOperator}[1]{%
  \mathord{\vcenter{\hbox{%
    \begin{tikzpicture}[
      x=6mm,y=5.5mm,draw=blue,line width=.7pt,
      line cap=round,line join=round]
      \foreach \x in {-3,-1,0,1,3}{
        \draw (\x,0)--(\x,2.8);
      }
      \foreach \y in {.25,2.45}{
        \foreach \x in {-2,2}{
          \node[font=\scriptsize,inner sep=0pt]
            at (\x,\y) {$\cdots$};
        }
      }
      \ifnum#1=0
        \draw[fill=white] (-3.35,.95) rectangle (.35,1.85);
        \node at (-1.5,1.4) {$f$};
      \else
        \draw[fill=white] (-.35,.95) rectangle (3.35,1.85);
        \node at (1.5,1.4) {$f$};
      \fi
      \foreach \x/\CouponStrandIndex in
        {-3/1,-1/{n-1},0/n,1/{n+1},3/{2n-1}}{
        \node[above,inner sep=2pt,font=\scriptsize]
          at (\x,2.8) {$\CouponStrandIndex$};
      }
    \end{tikzpicture}%
  }}}%
}
\newcommand{\LOperator}{\LabelledCouponOperator{0}}
\newcommand{\ROperator}{\LabelledCouponOperator{1}}

\begin{document}
\title{Property F and weak integrality of braided fusion categories}
\author{Enyu Li and Yuze Ruan}
\address{E. Li, Qiuzhen College, Tsinghua University, Beijing, 100084, China}
\email{ley24@mails.tsinghua.edu.cn}

\address{Y. Ruan, Yau Mathematical Sciences Center and Department of Mathematics, Tsinghua University, Beijing, 100084, China}
\email{yuzeruan@mail.tsinghua.edu.cn}

\begin{abstract}
We prove that a braided fusion category with property F is weakly integral, provided it admits a pivotal structure whose left categorical dimensions are positive on the intersection of its M\"uger center and its adjoint subcategory. In particular, the converse implication in
the property F conjecture holds for modular fusion categories and pseudounitary braided fusion categories. We also show that the positivity condition holds for every pivotal braided fusion category if and only if every braided fusion category over $\BC$ admits a pivotal structure.
\end{abstract}
\maketitle

\section{Introduction}
Braided fusion categories provide a natural source of braid group representations, with applications to quantum topology and topological quantum computation \cite{Tur94,BK01,RW18}. For an object $X$ in a braided fusion category $\CC$, the braiding defines a homomorphism
\[
\rho^X_{n}:B_n\longrightarrow \End_{\CC}(X^{\otimes n}).
\]
We say the object $X$ has property F if the image of $\rho^{X}_n$ is finite for all $n\geq 2$, The category $\CC$ is said to have property F if every object $X$ has property F. The property F conjecture, formulated in the modular setting in \cite{RSW09} and for braided fusion categories in \cite{NR11}, predicts that $\CC$ has property F if and only if it is weakly integral, that is, $\FPdim(\CC)\in\BZ$. It thus proposes an arithmetic characterization of a representation-theoretic property of the braiding.

Substantial progress has been made on the implication from weak integrality to property F. Finite-image results are known for braided group-theoretical categories \cite{ERW08} and for several families of metaplectic categories \cite{RW17,GRR20}. Green and Nikshych proved that every weakly group-theoretical braided fusion category has property F \cite{GN21}. In related work, the authors proved generalized property F for weakly group-theoretical modular categories with their canonical positive spherical structures, establishing finiteness of projective mapping class group images \cite{LiRuan2026}. These results concern the sufficiency of categorical arithmetic conditions for finiteness.

The converse direction has been verified in many quantum-group examples through the analysis of particular braid group representations \cite{NR11,Rowell2010}. Localization provides related arithmetic restrictions: if the braid representations associated with $X$ are localizable and their images span $\operatorname{End}_{\CC}(X^{\otimes n})$ for every $n$, then $\FPdim(X)^2\in\BZ$ \cite{RW2012}. Analogous restrictions have been established for generalized and quasi-localizations \cite{GHR2013}. In this paper, we give a uniform categorical argument for the necessary direction of the property F conjecture. Our main result applies to pivotal braided fusion categories whose pivotal dimensions are positive on the intersection of the M\"uger center and the adjoint subcategory.
\begin{Thm}
Let $\CC$ be a braided fusion category over $\BC$ equipped with a pivotal structure $j$. Suppose that
\begin{equation}\label{eq:key_assumption}
d_j^L(X)>0\qquad \text{for every nonzero }X\in \CC'\cap\CC_{\mathrm{ad}}.
\end{equation}
If $\CC$ has property F, then $\CC$ is weakly integral.   
\end{Thm}
Here $\CC'$ is the M\"uger center, and $\CC_{\mathrm{ad}}$ is the adjoint subcategory of $\CC$. Note that the positivity condition \eqref{eq:key_assumption} is independent of the chosen pivotal structure as all pivotal structures agree on $\mathcal C_{\mathrm{ad}}$. In particular, property F implies weak integrality for every modular fusion category and every pseudounitary braided fusion category. We also have
\begin{Cor}
If $\CC$ is a pseudounitary braided category and $X\in \Irr(\CC)$ has property F, then $\FPdim(X)^2\in\BZ$.   
\end{Cor}

The key step is the following trace estimate and the characterization of the equality: For a simple object $X$, put
\[
A_n=\mathrm{Span}_\BC\rho^X_{n}(B_n),\qquad \tau_n(T)=\frac{\operatorname{tr}_{X^{\otimes n}}(T)} {d(X)\FPdim(X)^{n-1}}.
\]Property F implies that the linear braid representations on the multiplicity spaces of tensor powers of $X$ are unitarizable for all $n$, choosing invariant inner product on these spaces makes $\{A_n\}_n$ a tower of finite dimensional $C^*$-algebras. We prove that
\[
|\tau_n(f)|\leq1
\]
for every unitary $f\in A_n$. More importantly, $|\tau_n|,n\ge 2$ attains 1 on unitary group of $A_n$ precisely when $|d(X)|=\FPdim(X)$ (Proposition \ref{prop:unitaryobj}). 

In a modular category, the simple object $X_{\mathrm{FP}}$ corresponding to the Frobenius--Perron character satisfies
\[
|\tau_2(c_{X_{\mathrm{FP}},X_{\mathrm{FP}}}^{\,2})|=1, \qquad d(X_{\mathrm{FP}})^2= \frac{\dim(\CC)}{\FPdim(\CC)}.
\]
Property F makes the relevant representations unitarizable, $u=c_{X_{\mathrm{FP}},X_{\mathrm{FP}}}^{\,2}$ attains the maximum, so the equality case implies that $\CC$ is pseudounitary. Applying the same argument to every Galois conjugate yields weak integrality.

For a general ribbon fusion category, we work in its center $\ZZ(\CC)$. The canonical braided tensor functor $\CC\boxtimes\CC^{\mathrm{rev}}\to\mathcal Z(\CC)$ has image $(\CC')'_{\mathcal Z(\CC)}$; when $\CC$ has property F, so does this image. The key step is to tensor $X_{\mathrm{FP}}\in\ZZ(\CC)$ with a suitable invertible object so that the resulting simple object lies in this image while preserving the required trace equality. We show that such an invertible object exists if and only if the positivity hypothesis \eqref{eq:key_assumption} holds (Corollary~\ref{cor:central_positivity_extension}). This equivalence follows from a pivotal extension criterion: a pivotal structure on a fusion subcategory $\DD\subseteq\CC$ extends to $\CC$ if and only if it agrees with the ambient pivotal structure on $\DD\cap\CC_{\mathrm{ad}}$ (Lemma~\ref{lem:pivotal_structure_extension}). Finally, we deduce the pivotal case by passing to the sphericalization.

We also relate the positivity hypothesis \eqref{eq:key_assumption} to the pivotal-structure existence conjecture \cite[Conj.~2.8]{ENO05}. We prove that this positivity condition holds for every pivotal braided fusion category if and only if every braided fusion category over $\BC$ admits a pivotal structure. Consequently, the braided case of the pivotal-structure existence conjecture would imply that every braided fusion category with property~F is weakly integral.

The paper is organized as follows. Section \ref{sec:prelim} recalls the necessary facts about fusion categories, braided and modular categories, Galois conjugation, and property F. Section \ref{sec:modular} establishes the trace estimate and treats the modular case. Section \ref{sec:braided} proves the pivotal extension criterion and the main theorem. Section \ref{sec:equivalence} discusses the equivalence with the pivotal-structure existence problem.

\section*{Acknowledgements}
YR is supported by a startup funding from YMSC, Tsinghua Shuimu Scholar Program and Beijing Municipal Science 
\& Technology Commission (Grant No. Z221100002722017).

The authors used generative AI to compute examples, examine patterns, and explore possible generalizations during the development of this work.  The authors take full responsibility for the mathematical contents of this paper.

\section{Preliminaries}\label{sec:prelim}
We recall the notation and standard facts used throughout the paper. For general background on fusion categories we refer to \cite{ENO05, EGNO}, and for braided and modular fusion categories to \cite{Tur94, BK01, DGNO}.
\subsection{Fusion category}
A \emph{fusion categories} over $\BC$ is a semisimple, $\BC$-linear abelian, rigid monoidal category with finite-dimensional Hom-spaces, finitely many isomorphism classes of simple objects, and a simple tensor unit $\one$.  For a fusion category $\CC$,  let $\Irr(\CC)$ denote a set of representatives of its simple objects, and write $X^*$ for the dual of $X$. Its Grothendieck ring is denoted by $K_0(\CC)$. The \emph{Frobenius–Perron dimension} $\FPdim(X)$ is the spectral radius of the fusion matrix of $X$. A \emph{pivotal structure} is a monoidal natural isomorphism $j:\operatorname{id}_{\CC}\xrightarrow{\sim}(-)^{**}$. Write $\tr^L_{j},\ \tr^R_j,\ d^L_j,\ d^R_j$ for the left and right traces and dimensions. The pivotal structure is \emph{spherical} if its left and right traces coincide. In that case, write $\tr_j,\ d_j(X)$ for the resulting spherical trace and dimension, omitting $j$ when the spherical structure is clear from the context. Then
\[
\dim(\CC)=\sum_{X\in\operatorname{Irr}(\CC)}d_j(X)^2,\qquad \FPdim(\CC)=\sum_{X\in\Irr(\CC)}\FPdim(X)^2
.\] 
Note that $\dim(\CC)$ is intrinsic. For every object $X$, both $d_j(X)$ and $\FPdim(X)$ are totally real algebraic integer. Both global dimensions are totally positive algebraic integers, and 
\[
|d_j(X)|\leq \FPdim(X),\qquad 1\leq \dim(\CC)\leq \FPdim(\CC).
\]
The category is \emph{pseudounitary} if $\dim(\CC)=\FPdim(\CC)$. Then it has a unique spherical structure $j_+$ with $d_{j_+}(X)=\FPdim(X)>0$. The category is called \emph{weakly integral} if $\FPdim(\CC)\in\BZ$. Every weakly integral fusion category is pseudounitary, and its simple objects satisfy $\FPdim(X)^2\in\BZ$; see \cite[Sections~8.3--8.5]{ENO05}.

The \emph{adjoint subcategory} $\CC_{\mathrm{ad}}$ is the fusion subcategory generated by the simple constituents of $X\otimes X^*$, where $X\in\Irr(\CC)$. The \emph{pointed subcategory} $\CC_{\mathrm{pt}}$ is the fusion subcategory generated by the invertible objects of $\CC$. The universal grading of the Grothendieck ring $K_0(\CC)$ induces a canonical faithful grading \cite{GN08}:
\[
\CC=\bigoplus_{g\in U(\CC)}\CC_g,
\qquad
\CC_e=\CC_{\mathrm{ad}},
\]
where $U(\CC)$ is called the \emph{universal grading group}. A nonzero object $X\in\CC_g$ is homogeneous of degree $\deg(X)=g$; every simple object is homogeneous. For nonzero homogeneous $X,Y$,
\[
 \deg(X\otimes Y)=\deg(X)\deg(Y),\qquad \deg(X^*)=\deg(X)^{-1},\qquad \deg(\one)=e.
\]
Every faithful group grading of $\CC$ is obtained from this grading via a unique surjective homomorphism from $U(\CC)$. There is a canonical identification
\begin{equation}\label{eq:piv_character}
\Aut_{\otimes}(\Id)\simeq \Hom(U(\CC),\BC^{\times}).
\end{equation}
If $\CC$ admits a pivotal (respectively, spherical) structure, then the set of such structures form a torsor over $\Hom(U(C),\BC^{\times})$ (respectively, $\Hom(U(C),\pm1)$), the character $\chi$ modifies
dimensions by $d^L_{j^\chi}(X)=\chi(\deg X)d^L_j(X)$ for homogeneous $X$. In particular, all pivotal structures agree on $\CC_{\mathrm{ad}}$.

Every fusion category $\CC$ admits a canonical \emph{sphericalization}. Its objects are pairs $(X,\varphi_X)$, where $X\in\Obj(\CC)$ and $\varphi_X:X\to X^{**}$ is an isomorphism satisfying $(\varphi_X)^{**}\circ\varphi_X=\gamma_X$, with $\gamma:\Id_{\CC}\to(-)^{****}$ the canonical Radford isomorphism. A morphism $a:(X,\varphi_X)\to(Y,\varphi_Y)$ is a morphism $a:X\to Y$ in $\CC$ satisfying $\varphi_Y\circ a=a^{**}\circ\varphi_X$. The maps $\varphi_X$ define a canonical spherical structure on $\widetilde{\CC}$. The forgetful tensor functor
$F:\widetilde{\CC}\to\CC$ is faithful and essentially
surjective, and sends simple objects to simple objects. Every simple $X$ has two simple lifts, represented by $(X,\varphi_X)$ and $(X,-\varphi_X)$, and
\[
\dim(\widetilde{\CC})=2\dim(\CC),\qquad
\FPdim(\widetilde{\CC})
=2\FPdim(\CC).
\]
If $\CC$ is braided, its braiding lifts to $\widetilde{\CC}$, making $F$ braided; also see \cite[Sec.~5.1]{Ostrik15}

\subsection{Braided and modular fusion categories}
A \emph{braiding} on $\CC$ is a natural isomorphism $c_{X, Y}: X \otimes Y \xrightarrow{\cong} Y \otimes X$ satisfying the hexagon axioms. A fusion category equipped with a braiding is called a \emph{braided fusion category}. A \emph{twist} is a natural isomorphism $\theta: \Id_{\CC}\to \Id_{\CC}$ satisfying the \emph{balancing equation}
\begin{equation}\label{eq:Balance_equ}
\theta_{X\otimes Y}=(\theta_X\otimes \theta_Y)\circ c_{Y,X}\circ c_{X,Y},\qquad \theta_{\one}=\Id_{\one}.
\end{equation}
It is a \emph{ribbon structure} if $\theta_{X^*}=(\theta_X)^*$. For $X\in\Irr(\CC)$, we also write $\theta_X$ for its twist scalar. These scalars are roots of unity by Vafa's theorem \cite{Etingof2002Vafa}. On a braided fusion category, spherical structures correspond to ribbon structures through the Drinfeld isomorphism $u_X:X\to X^{**}$, and $j_X=u_X\theta_X$. A braided fusion category equipped with ribbon or spherical structure is called a \emph{ribbon fusion category}. 

The \emph{Müger centralizer} of a fusion subcategory $\DD\subseteq\CC$ is the fusion subcategory
\[
\DD_{\CC}'=\left\{X\in\CC\ \middle|\ c_{Y,X}\circ c_{X,Y}=\operatorname{id}_{X\otimes Y}\text{ for every }Y\in\DD\right\}.
\]
A ribbon fusion category $\CC$ is \emph{modular} if $\CC'_{\CC}=\operatorname{Vec}$, equivalently if its unnormalized S-matrix is invertible. We use
$$
S_{X,Y} := \Tr_{X\otimes Y}(c_{Y, X}\circ c_{X, Y}),\qquad T_{X,Y} := \delta_{X, Y}\cdot\theta_X ,\qquad X,Y\in\Irr(\CC).
$$
In particular, $S_{X, \one} = S_{\one, X} = d_j(X)$ and $T$ has finite order. By the Verlinde formula, the fusion matrices are simultaneously diagonalized by S-matrix. For each $X\in \Irr(\CC)$, we have the correspooding character:
\[
\lambda_X: K_0(\CC)\longrightarrow\BC, \qquad \lambda_X(Y)=\frac{S_{XY}}{d(X)},
\]
and $\{\lambda_X\mid X\in\Irr(\CC)\}$ forms the complete set of distinct unital ring characters. Consequently, there exists a unique self-dual simple object $X_{\mathrm{FP}}$ satisfying
\begin{equation}\label{eq:X_FP}
 \lambda_{X_{\mathrm{FP}}}(Y)=\FPdim(Y), \qquad d_j(X_{\mathrm{FP}})^2=\frac{\dim(\CC)}{\FPdim(\CC)}; 
\end{equation}
compare to \cite[Prop.~8.22]{ENO05}.

In the modular case, group in \eqref{eq:piv_character} is canonically isomorphic to $\Inv(\CC)$. Under this isomorphism, $v\in\Inv(\CC)$ corresponds to the tensor natural automorphism $\eta^v$ of $\id_{\CC}$ with components (\cite[Lem.~8.22.9]{EGNO}),
\begin{equation}\label{eq:inv_aut}
(\eta^v)_Y=\frac{S_{v,Y}}{d(v)d(Y)}\,\id_Y, \qquad Y\in\Irr(\CC).
\end{equation}

The \emph{Drinfeld center} $\ZZ(\AA)$ of a fusion category $\AA$ consists of pairs $(X,\gamma)$, where $\gamma_Y:X\otimes Y\xrightarrow{\sim}Y\otimes X$ is a half-braiding. If $\AA$ is spherical, then $\ZZ(\AA)$ is modular with the induced spherical structure, for which the forgetful functor is pivotal. In particular,
\[
d_{\ZZ(\AA)}(X,\gamma)=d_{\AA}(X),\qquad
\dim(\ZZ(\AA))=\dim(\AA)^2,
\]
and $\FPdim(\ZZ(\AA))
=\FPdim(\AA)^2$.
% There is a canonical group injection
% \[
% \Hom(U(\CC),\BC^\times)\simeq\Aut_{\otimes}(\id_{\CC}) \hookrightarrow\Inv(\ZZ(\CC)),\qquad
% \chi\longmapsto[(\one,\gamma^\chi)].
% \]
% Here $\gamma^\chi_Y=\chi(\deg(Y))\Id_Y$ after the unit identifications; its image consists precisely of the invertible central objects whose underlying object is isomorphic to $\one$ 
For ribbon $\CC$, equip $\CC^{\rev}$ with $c^{\rev}_{X,Y}=c_{Y,X}^{-1}$ and $\theta_X^{\rev}=\theta_X^{-1}$.
With the induced ribbon structure on $\ZZ(\CC)$, the canonical embeddings are
\[
\iota_+:\CC\hookrightarrow\ZZ(\CC),\quad X\mapsto(X,c_{X,-}), \qquad
\iota_-:\CC^{\rev}\hookrightarrow\ZZ(\CC),\quad X\mapsto(X,c_{-,X}^{-1})
\]
Their images centralize one another and induce a ribbon functor
$G:\CC\boxtimes\CC^{rev}\to\ZZ(\CC)$. Identify the M\"uger center $\CC'=\CC'_{\CC}$ with its common image $\iota_+(\CC')=\iota_-(\CC')$ in $\ZZ(\CC)$. Then 
\[
 \operatorname{Im}(G)=(\CC')'_{\ZZ(\CC)}.
\]

\subsection{Galois conjugation and property F}
For $\sigma\in \Aut(\BC)$, a \emph{Galois conjugate} $\CC^{\sigma}$ is obtained by applying $\sigma$ to all structural coefficients of $\CC$, including the chosen pivotal, braiding, and ribbon data. Fusion rules, Frobenius-Perron dimensions, and the universal grading are unchanged, whereas
\begin{equation}
d_{j^{\sigma}}(X^{\sigma})=\sigma(d_j(X)),\quad \dim(\CC^{\sigma})=\sigma(\dim(\CC)).   
\end{equation}

For $X\in\Obj(\CC)$, the braiding gives representations
\[
\rho_n^{X}:B_n\to \End_{\CC}(X^{\otimes n}), \qquad \sigma_i\mapsto\Id_{i-1}\otimes c_{X,X}\otimes \Id_{n-i-1}. 
\]
Semisimplicity gives the natural isomorphism
\[\End_{\CC}(X^{\otimes n})\cong\bigoplus_{Y\in\Irr(\CC)}{\End}_{\BC}\left({\Hom}_{\CC}(Y,X^{\otimes n})\right)\]
 gives linear representations \[\rho_{n,Y}^{X}:B_n\to{\End}_{\BC}\left({\Hom}_{\CC}(Y,X^{\otimes n})\right).\]

An object has \emph{property F} when these images are finite for every $n$; the category has property F when every object does \cite[Definition~2.1]{NR11}. Property F passes to subobjects and fusion subcategories, and is preserved by reversing the braiding and by Deligne tensor products. It also passes to the image through braided tensor functors. If a simple object $X$ has property F, then so does the fusion subcategory generated by it; compare \cite[Section~2]{NR11}. Finally, applying $\sigma$ entrywise identifies a finite braid image with
its Galois-conjugate image. Thus $X$ has property F if and only if $X^\sigma$ does, and $\CC$ has property F if and only if $\mathcal C^\sigma$ does. Property F depends only on the braiding, and is independent of the choice of pivotal or spherical structure.

\section{Trace estimates and the modular case}\label{sec:modular} In this section we shall prove that a modular category $\CC$ with property F is weakly integral (see theorem \ref{thm:modular_case}). The idea is that a representation with finite image is always unitarizable, which puts restrictions on eigenvalues, in particular, a bound to the trace of images of linear representations of braids. By a careful analysis of the tracial functional, we shall show that the trace bound on images of $\rho^X_n(n\ge 2)$ is attained precisely when $\FPdim(X)=|d(X)|$ (see proposition \ref{prop:unitaryobj}). For a modular category, the trace bound is attained at $c_{X_{\mathrm{FP}},X_{\mathrm{FP}}}^2$ where $X_{\mathrm{FP}}$ is the simple object corresponding to the Frobenius–Perron character. It follows that $|d(X_{\mathrm{FP}})|=\FPdim(X_{\mathrm{FP}})$, which implies that every modular category with property F is pseudounitary. Weak integrality then follows immediately from the following lemma:  
\begin{lem}\label{galois} A fusion category $\CC$ is weakly integral if and only if every Galois conjugate of $\CC$ is pseudounitary.
\end{lem}\begin{proof} Let \(\CC^{\sigma}\) be the category obtained by conjugate the structural coefficients of $\CC$ using  \(\sigma\in\mathrm{Aut}({\BC})\). 
\[\sigma(\dim(\CC))=\dim(\CC^{\sigma})=\FPdim (\CC^{\sigma})=\FPdim{(\CC)}=\dim(\CC)\]
    which implies that $\dim (\CC)$ is a rational algebraic integer, hence an integer. The converse is trivial. 
\end{proof}

Property F implies that the linear representations of braid factors through a finite group, such a representation is always unitarizable.   

\begin{defn}
We say $X$ is {unitarizable} if there is a inner product on ${\Hom}_{\CC}(Y,X^{\otimes n})$ for every $n$ and $Y\in\Irr(\CC)$ such that the associated braid group representations \[\rho_{n,Y}^{X}:B_n\to{\End}_{\BC}\left({\Hom}_{\CC}(Y,X^{\otimes n})\right)\] is unitary for every $n$ and $Y\in \Irr(\CC)$.
\end{defn} If the category is unitary, then every object is unitarizable, but not the converse.  Also note that our definition is slightly different from that of \cite[Def.~2.10]{Rowell2005}. When $X$ is unitarizable, there is a $C^*$-algebra structure on algebra spanned by braid images so that tensor products preserve unitarity:

\begin{lem}\label{lem:staralg}
        Let $X$ be a simple unitarizable object in a ribbon fusion category. Let $A_n=\mathrm{Span}_\BC(\rho_n^X(B_n))\subset \End_{\CC}(X^{\otimes n})$. Then there is a $C^*$-algebra structure on $A_n$ for every $n$ such that the tensor product $A_k\otimes A_l\to A_{k+l}$ is a $*$-homomorphism and  preserves unitary elements.
\end{lem}

\begin{proof}
    Since $X$ is unitarizable, there is a $B_n$ invariant inner product on $\Hom_\CC(Y,X^{\otimes n})$, which makes $\End_\BC(\Hom_\CC(Y,X^{\otimes n}))$, hence $\End_\CC(X^{\otimes n})$, a $C^*$-algebra and $A_n$ a $C^*$-subalgebra. Note that $\rho_n^X(b)^*=\rho_n(b^{-1})$ for $b\in B_n$, $(\rho_k^X(b_1)\otimes\rho_l^X(b_2) )^*=\rho_{k+l}^X(b_1\otimes b_2)^*=\rho_{k+l}^{X}(b_1^{-1}\otimes b_2^{-1})=\rho_{k}^X(b_1)^*\otimes \rho_l^X(b_2)^*$, hence $A_k\otimes A_l\to A_{k+l}$ is a $*$-homomorphism and hence preserves unitary elements.
\end{proof}

\begin{rmk}
It is crucial that we restrict ourselves to the $*$-subalgebra generated by braids, the tensor product  $\End(X^{\otimes k})\otimes  \End(X^{\otimes {l}})\to \End(X^{\otimes (k+l)})$  may not be $*$-homomorphism.      
\end{rmk}

The following characterization of the maximal loci of a trace functional over a $C^*$-algebra is crucial to our proof of proposition \ref{prop:unitaryobj}:

\begin{lem}\label{lem:maxtrace}
    Let $A$ be a finite dimensional $C^*$-algebra. Let $t:A\to \BC$ be a tracial functional, and \[\mathrm{Rad}_t=\{x\in A: t(xy)=0 \text{ for every }  y\in A\}.\] If $|t|$ attains its maximum on unitary subgroup at $u$, then $u+\mathrm{Rad}_t\in Z(A/\mathrm{Rad}_t)$.
\end{lem}
\begin{proof}Write \[A\cong\bigoplus_\alpha \mathrm{Mat}_{n_\alpha}{\BC},\qquad t=\sum_\alpha w_\alpha\Tr_\alpha,\]then\[\mathrm{Rad}_t\cong\bigoplus_{w_\alpha=0}\mathrm{Mat}_{n_\alpha}(\mathbb C), x\mapsto (x_{\alpha})_{\alpha}
\]Note that for unitary $x$:
\[|t(x)|=|\sum_{\alpha}w_\alpha \Tr_\alpha(x_{\alpha})|\leq \sum_\alpha |w_\alpha|n_\alpha.\] with equality if and only if $x_\alpha= w|w_\alpha|/w_\alpha I_{n_\alpha}$ for $w_\alpha\neq 0$ and $|w|=1$. It's then clear that $x$ is central modulo $\mathrm{Rad}_t$.\end{proof}

Fix a simple object $X$. For convenience we renormalize the trace functional on $\End_\CC(X^{\otimes n})$ by \[
\tau_n=\frac{\tr_{X^{\otimes n}}}{d(X)\FPdim(X)^{n-1}}
\]
so that for $S\in \End{(X^{\otimes k})}, T\in \End({X^{\otimes l}})$
\[\tau_{k+l}(S\otimes T)=\frac{d(X)}{\FPdim(X)}\tau_k(S)\tau_l(T)\] and\[
\tau_{k+l-1}(S\#T)=\tau_{k}(S)\tau_l(T).
\] where 
\[S\# T:=(S\otimes \id_{X^{\otimes (l-1)}})\circ(\id_{X^{\otimes (k-1)}}\otimes T )=\OverlapCoupons{S}{T}.\]

Under the identification\[{\End}_{\CC}(X^{\otimes n})
\xrightarrow{\sim}
\bigoplus_{Y\in\operatorname{Irr}(\CC)}
{\End}_{\BC}\!\left({\Hom}_{\CC}(Y,X^{\otimes n})
\right), T \mapsto \oplus_Y T_Y,
\] where $T_Y$ is the linear operator induced by postcomposition with \(T\), the trace functional is
\[
\begin{aligned}
  {\tr}(T)
  &=\TraceDiagram{0}   =     \sum_{Y,a}\TraceDiagram{2}\\&=  \sum_{Y,a}(T_Y)_{aa}\,\TraceDiagram{3}
   \\&=\sum_Y d(Y)\,\operatorname{Tr}_{Y}(T_Y).
\end{aligned}
\]where $\Tr_Y$ is the matrix trace. Thus \[\tau_n(T)=\frac{1}{d(X)\FPdim(X)^{n-1}}\sum_Y d(Y)\Tr_Y(T_Y).\]

\begin{lem}\label{braketX}
    Let $X$ be simple in a spherical fusion category. If $\FPdim(X)=\epsilon d(X),\epsilon=\pm1$, then $\epsilon^{n+1}\tau_n$ on $\End_\CC(X^{\otimes n})$ is positive and the subcategory  $\langle X\rangle$ generated by $X$ is pseudounitary. If $Y$ is simple subobject of $X^{\otimes n}$, $d(Y)=\epsilon^n \FPdim(Y)$.  
\end{lem}
\begin{proof}
    Suppose $d(X)=\epsilon  \FPdim(X)$ for $\epsilon=\pm 1$. If $X^{\otimes n}=\bigoplus_{Y\in \Irr(\CC)} Y^{\oplus N_Y}, N_Y>0$,  \[  0=\FPdim(X)^{n}-\epsilon^{n} d(X)^{n}=\sum_Y (\FPdim(Y)-\epsilon^{n} d(Y))N_Y \]thus $\FPdim(Y)= \epsilon^n d(Y)$ and for $T\in\End_\CC(X^{\otimes {n}})$, \[\epsilon^{n+1}\tau_{n}(T)=\frac{\epsilon^{n+1}}{d(X)\FPdim(X)^{n-1}}\sum_{Y}d(Y)\Tr_Y(T_Y)=\frac{1}{\FPdim(X)^{n}}\sum_{Y}\FPdim(Y)\Tr_Y(T_Y)\] which is positive. 
\end{proof}

Now we can prove the following core properties of the trace $\tau_n$ on $A_n$:

\begin{prop}\label{prop:unitaryobj}
    Let $X$ be a simple unitarizable object in a ribbon fusion category. Let $A_n=\mathrm{Span}_\BC(\rho_n^X(B_n))$. Then for unitary $f\in A_n$: \[|\tau_n(f)|\le 1.\] If $n\geq 2$, the bound is attained by some unitary \(f\in A_n\) if and only if $d(X)=\pm \FPdim(X)$. In that case, $f$ is scalar. 
\end{prop}

\begin{proof}
Since $X$ is unitarizable, $T_Y$ has eigenvalue with modulus one for every unitary $T\in A_n$, in particular we have the trace bound:
\begin{align*}
|\tr(T)|=&\left|\sum_Y d_Y\Tr(T_Y)\right|\leq\sum_{Y\in{\Irr(\CC)}} |d(Y)|\dim{\Hom}_{\CC}(Y,X^{\otimes n})\\\leq &\sum_{Y\in{\Irr(\CC)}} \FPdim(Y)\dim{\Hom}_{\CC}(Y,X^{\otimes n})\\=& \FPdim(X)^{n}
\end{align*}
thus \[|\tau_n(T)|\leq \frac{\FPdim(X)}{|d(X)|}\] i.e. our $\tau_n$ is uniformly bounded on unitary elements for every $n$.
By lemma \ref{lem:staralg}, $T^{\# m}$ is also unitary in $A_{nm-m+1}$, \[|\tau_{nm-m+1}(T^{\#m})|=|\tau_{n}(T)|^m\le \frac{\FPdim(X)}{|d(X)|}.\] This forces $|\tau_n(T)|\leq 1$ for every unitary $T\in A_n$.

% The dimensions $d(X),d(Y)$ are real, which implies $\tau$ is Hermitian.
Now suppose $n\ge 2$ and $|\tau_n(f)|=1$ for some $n$. Note that $|\tau_{2n-1}(f\# f)|=|\tau_n(f)|^2=1$, so $|\tau_{2n-1}|$ attains its maximum at $f\# f$ on the unitary elements of $A_{2n-1}$. For convenience write \[L:=\LOperator
\qquad
R:=\ROperator.\] By lemma \ref{lem:maxtrace}, $f\# f=LR$ is central modulo $\mathrm{Rad}_{\tau_{2n-1}}$ in $A_{2n-1}$. If $n\ge 3$, the braid generators $\sigma_j, j\neq n-1,n$ commutes with $L$ or $R$, thus they commute with both modulo $\mathrm{Rad}_{\tau_{2n-1}}$. It follows that
\[\CommuteDiagram{0}
= \CommuteDiagram{1}
= \CommuteDiagram{2}
= \CommuteDiagram{3}
  \quad\mod \mathrm{Rad_{\tau_{2n-1}}}\] and
\[\ForwardCommuteDiagram{0}
= \ForwardCommuteDiagram{1}
= \ForwardCommuteDiagram{2}
= \ForwardCommuteDiagram{3}
  \quad\mod \mathrm{Rad_{\tau_{2n-1}}}.\]
  As a consequence, $\sigma_{n-1},\sigma_n$ also commute with $L$ and $R$, hence both $L$ and $R$ are central modulo $\mathrm{Rad}_{\tau_{2n-1}}$. If $n=2$, braid generators commute with $L$ or $R$ since $A_2$ is commutative, again using $LR$ is central modulo $\mathrm{Rad}_{\tau_{2n-1}}$, they must commute with both modulo $\mathrm{Rad}_{\tau_{2n-1}}$.
  We have:
\[
  \CouponSlide{1}
 \quad {=}\quad
  \CouponSlide{2}
  \quad {=}\quad
  \CouponSlide{3}
  \quad\mod \mathrm{Rad_{\tau_{2n-1}}}\]
  so
  \begin{equation}
    \label{eq:LR}
L=  \CouponOp{0}=\CouponOp{1}=R \mod \mathrm{Rad_{\tau_{2n-1}}},
\end{equation}
as a consequence 
\[f\# f=LR=L^2    \mod \mathrm{Rad_{\tau_{2n-1}}}\]
and 
\[1=|\tau_{2n-1}(f\# f)|=|\tau_{2n-1}(f^2\otimes \id_{X^{\otimes (n-1)}})|=\frac{|d(X)|^{n-1}}{\FPdim(X)^{n-1}}|\tau_n(f^2)|\le\frac{|d(X)|^{n-1}}{\FPdim(X)^{n-1}}\le 1.
\] Hence $|d(X)|=\FPdim(X)$.

Suppose $d(X)=\epsilon\FPdim(X)$ for $\epsilon=\pm 1$, $\tau_{2n-1}$ is positive by lemma \ref{braketX}. It follows that the restriction of $\tau_{2n-1}$ on $A_{2n-1}$ is also positive, in particular $\mathrm{Rad}_{\tau_{2n-1}}=0$. By equation \eqref{eq:LR}, $L=R$, closing the last $n-1$ strands gives \[f=\frac{\tr(f)}{d(X)^n}\id_{X^{\otimes n}}. \]
% Since 
% \[\sum_Y N_Y\FPdim(Y)=\FPdim(X)^{n}=|\tr(f)|=\left|\sum_{Y}\Tr(T_Y)d(Y)\right|=\left|\sum_Y\Tr(T_Y)\FPdim(Y)\right|.\] Since $T_Y$ is unitary, the equality holds only if $T_Y=\lambda\id_{Y^{\otimes N_Y}}$ for some $|\lambda|=1$. It follows that $T=\lambda \id$.
Thus $|\tau_n|$ attains $1$ precisely at scalars with modulus one.
\end{proof}
\begin{rmk}
The trace is always nondegenerate on $\End_\CC(X^{\otimes {2n-1}})$, but its restriction on $A_{2n-1}$ might have nonzero radical. If $\FPdim(X)=|d(X)|$, $\tau_{2n-1}$ is positive, hence the restriction on $A_{2n-1}$ has no radical. 
\end{rmk}  
% In this case we have a faithful positive tracial functional on $A_n$. If $A_n$ is semisimple, the braid generators $\sigma_i$ gives an injective map $B_n\to \BC^{n},\beta\mapsto \tau(\beta\sigma_i)$. Thus if we know $\tau$ is bounded on braid images, $B_n$ has compact image hence must be unitary.  
\begin{prop}\label{prop:tracebound}
    Let $X$ be a simple object in a ribbon fusion category. If there is some $f\in \rho_n^X(B_n), n\ge 2$ such that 
    $$|\tau_n(f)|=1.$$
Consider the following conditions:
\begin{enumerate}
  \item $X$ has property~$F$;
  \item $X$ is unitarizable;
  \item $f$ is scalar;
  \item $\FPdim(X)=|d(X)|$.
\end{enumerate}
Then
$\textup{(1)}\Rightarrow\textup{(2)} \Rightarrow\textup{(3)}\Rightarrow\textup{(4)}$. If $n=2$ and $f=c_{X,X}^{2}$, all four conditions are equivalent.\end{prop}
\begin{proof}
(1)$\Rightarrow$(2): For any inner product $\langle-,-\rangle$ on ${\Hom}_{\CC}(Y,X^{\otimes n}),$ \[(x,y):=\sum_{g\in \rho^X_{n,Y}(B_n)}\langle gx,gy\rangle\] is a $B_n$-invariant inner product.

(2)$\Rightarrow$(3) and (4) is proposition \ref{prop:unitaryobj}.  

(3)$\Rightarrow$(4): By Vafa's theorem \cite[Thm.4.1]{Etingof2002Vafa}, every braid generator has finite order, hence $\det (f|_{\Hom_\CC(Y,X^{\otimes n})})$ has finite order. It follows that $f$ is a root of unity and $|\tr(f)|=|d(X)|^{n}=|d(X)|\FPdim(X)^{n-1}$. Hence $|d(X)|= \FPdim(X)$. 

(4)$\Rightarrow$(3)$\Rightarrow$(1) When $f=c_{X,X}^2$: if $X^{\otimes 2}=\bigoplus_Y Y^{\oplus N_Y}$, then $d(Y)=\FPdim(Y)$ by lemma \ref{braketX}, so 
\begin{align*}
    \left|\sum_Y {\FPdim(Y)N_Y} \right|=&\FPdim(X)^2=|\tr(f)|\\=&\left|\sum_Y {d(Y)}\Tr_Y((c_{X,X}^2)_Y) \right|\\=&\left|\sum_Y {\FPdim(Y)}\Tr_Y((c_{X,X}^2)_Y) \right|\end{align*}
    Each $(c_{X,X}^2)_Y$ is a root of unity by Vafa's theorem, thus the equality holds only if 
 $c_{X,X}^2=\lambda\id_{X^{\otimes 2}}$ with $\lambda$ root of unity. If $\sigma_i\in B_n$ are braid generators, $\hat\rho(\sigma_i)=\rho^{X}_n(\sigma_i)/\sqrt\lambda$ defines a braid representation with $\left(\hat\rho(\sigma_i)\right)^2=\id_{X^{\otimes n}}$, thus $\hat\rho$ factor through the symmetric group $\mathrm{Sym}_n$, hence has finite image. It follows that $\rho_n^X(B_n)\subset \hat\rho(B_n)\langle \sqrt\lambda\id_{X^{\otimes n}}\rangle $ and hence $\rho_n^X$ has finite image.  
\end{proof}

\begin{Thm}\label{thm:modular_case}
Every modular category $\CC$ with property F is pseudounitary, hence weakly integral. 
\end{Thm}
\begin{proof} By (\ref{eq:X_FP}), a modular category has a simple object $X_\mathrm{FP}$ satisfying 
\begin{equation}\label{eq:key_condition}
|\tr(c_{X,X}^2)|=|d(X)|\FPdim(X), d(X)^2=\frac{\dim(\CC)}{\FPdim(\CC)}.   
\end{equation}
Proposition \ref{prop:tracebound} implies 
\[1\leq \FPdim(X)^2= d(X)^2 =\frac{\dim(\CC)}{\FPdim(\CC)}\leq 1\]
thus $\CC$ is pseudounitary. Since every Galois conjugate $\CC^{\sigma}$ is modular with property F, lemma \ref{galois} then implies that $\CC$ is weakly integral.
\end{proof}

\begin{rmk}
    The result of proposition \ref{prop:unitaryobj} and \ref{prop:tracebound} holds even if $\CC$ has no pivotal structure, with $d(X)$ replaced by $\delta(X)=\sqrt{|X|^2}=\sqrt{d^L(X)d^{R}(X)}$ and the normalized trace replaced by 
\[\tau_{n}=\frac{\delta(X)^{n-1} }{\FPdim(X)^{n-1}d^L(X)^{n}}{\tr}^L=\frac{\delta(X)^{n-1}}{\FPdim(X)^{n-1}d^R(X)^{n}}  {\tr}^R\]which is independent of the chosen isomorphism $X\to X^{**}$ to define $\tr^L,\tr^R$ and $d^L,d^R$. There is also an intrinsic $S$ matrix (See \cite[Sec.~8.19]{EGNO}).  
\[
S_{X,Y}=\frac{\delta(X)\delta(Y)}{d^R(X)d^L(Y)}\left({\tr}^R\otimes{\tr}^L\right)\left(c_{Y,X}\circ c_{X,Y}\right).\]
However, there might be no simple object $X$ satisfying \[\frac{|S_{X,Y}|}{\delta(X)}={\FPdim(Y)}\]for every simple $Y$. In fact, in a nondegenerate braided fusion category, such an object exists if and only if $\CC$ admits a pivotal structure.
\end{rmk}

As illustrated in Theorem \ref{thm:modular_case}, it is crucial to find an object $X$ satisfying the condition \eqref{eq:key_condition}. For this purpose, the object $X_{\mathrm{FP}}$ is not the unique one.

\begin{lem}\label{lem:perturb}
Let $\CC$ be a modular fusion category, The condition \eqref{eq:key_condition} holds for $X_{\mathrm{FP}}\otimes v$ for any $v\in \Inv(\CC)$.    
\end{lem}
\begin{proof}
Consider the object $Z_v=:X_{\mathrm{FP}}\otimes v$ for any $v\in \Inv(\CC)$, then $|d(Z)|=|d(X_{\mathrm{FP}})|$ we have
\[
S_{Z,Y}=\frac{S_{X_{\mathrm{FP}},Y}S_{v,Y}}{d(Y)}=d(X_{\mathrm{FP}})\FPdim(Y)\frac{S_{vY}}{d(Y)},
\]
and $\frac{S_{vY}}{d(Y)}$ is a root of unity as $v$ is invertible. Therefore, taking the absolute values gives the desired result.
\end{proof}

\section{Pivotal extensions and the braided case}\label{sec:braided}
In this section, we will generalize the result to the braided fusion category, under the condition that $\CC$ has a pivotal structure which is positive restricted to $\CC'\cap\CC_{\mathrm{ad}}$, and note that the positivity condition is independent of the pivotal structure. The idea of the proof is as follows:
\begin{itemize}
\item Following the discussion in the preliminaries, the dominant ribbon functor $G: \CC\boxtimes \CC^{rev}\to i_+(\CC)\vee i_-(\CC)=(\CC')'_{\ZZ(\CC)}$ transforms the property F of $\CC$ to $(\CC')'_{\ZZ(\CC)}$, a subcategory of a modular category. 
\item From the proof of Theorem \ref{thm:modular_case} in the modular case, it suffices to find objects in $(\CC')'_{\ZZ(\CC)}$ satisfying condition \eqref{eq:key_condition} for every Galois conjugation, indeed it will imply $\ZZ(\CC)$ is weakly integral, thus so is $\CC$. therefore the idea is to perturb $X_{\mathrm{FP}}$ by an invertible object as in Lemma \ref{lem:perturb} such that the new object lands in the image of $G$. This step will have obstructions described in the following Lemmas.
\end{itemize}

\begin{lem}\label{lem:perturb_and_extension}
Let $\CC$ be a modular fusion category with corresponding spherical structure $j$, $\DD\subseteq\CC$ a pseudounitary fusion subcategory and $Z_v=X_{\mathrm{FP}}\otimes v$ for $v\in \Inv(\CC)$ as described in Lemma \ref{lem:perturb}, the following two statements are equivalent:
\begin{enumerate}
    \item There exists $v\in\Inv(\CC)$, such that $Z_{v}\in\DD'$, 
    \item The positive spherical structure $j_+$ on $\DD$ extends to a pivotal structure on $\CC$.
\end{enumerate} 
\end{lem}

\begin{proof}
As computed in Lemma \ref{lem:perturb}, we have \begin{align*}
Z_v=X_\mathrm{FP}\otimes v\in (\DD)'_{\CC}&\Leftrightarrow S_{Z_v,Y}=d_j(Z_v)d_j(Y),\quad Y\in\Irr(\DD), \\
&\Leftrightarrow d_j(Y)=d_{j_+}(Y)\frac{S_{vY}}{d(v)d(Y)},\quad Y\in\Irr(\DD).\\
&\xLeftrightarrow{\eqref{eq:inv_aut}} j\circ (\eta^v)^{-1}|_{\DD}=j_+,
\end{align*}
where the first equivalence follows from \cite[Prop.~C.3]{DGNO}.
\end{proof}

Let $\CC$ be a ribbon fusion category, by Deligne's theorem \cite{Del02}, $\CC'$ is integral and thus pseudounitary. Identify $\CC'$ with the its image in $\ZZ(\CC)$. Apply Lemma \ref{lem:perturb_and_extension} to $\ZZ(\CC)$ and its pseudounitary subcategory $\CC'$, we have

\begin{Cor}\label{cor:perturb_and_extension_center}
The following two statements are equivalent:
\begin{enumerate}
    \item There exists $v\in\Inv(\ZZ(\CC))$, such that $Z_{v}\in(\CC')'_{\ZZ(\CC)}$, $(Z_v=X_{\mathrm{FP}}\otimes v)$
    \item The positive spherical structure $j_+$ on $\CC'$ extends to a pivotal structure on $\CC$.
\end{enumerate}     
\end{Cor}
\begin{proof}
Condition (1) is equivalent to the existence of a pivotal structure on $\ZZ(\CC)$ extending $j_+$, which restricts along $\iota_+:\CC\hookrightarrow\ZZ(\CC)$ to a pivotal structure on $\CC$ extending $j_+$, Conversely, any pivotal extension to $\CC$ induces a pivotal structure on $\ZZ(\CC)$ with the same restriction to $\CC'$.    
\end{proof}

We next give a criterion for extending a pivotal structure from a fusion subcategory to the ambient category. Let $\CC$ be a braided fusion category, equipped with a pivotal structure $j$, and let $\DD\subseteq\CC$ be a fusion subcategory equipped with a pivotal structure $j_{\DD}$. Write their universal grading as
\[
\CC=\bigoplus_{g\in U(\CC)}\CC_g, \qquad \DD=\bigoplus_{g\in U(\DD)}\DD_g.
\]

The universal property of the grading of $\DD$ gives a natural homomorphism
\[
\phi:U(\DD)\longrightarrow U(\CC),\qquad \DD_u\subseteq\CC_{\phi(u)}.
\]
Set
\begin{align*}
H&:=\operatorname{im}(\phi)=\{h\in U(\CC):\DD\cap\CC_h\ne0\},\\
K&:=\ker(\phi)=\{u\in U(\DD):\DD_u\cap\CC_e\ne0\}.
\end{align*}
Here a nonzero intersection means that it contains a nonzero object. There is a short exact sequence
\[
1\longrightarrow K\longrightarrow U(\DD)\xrightarrow{\phi}H\longrightarrow1.
\]
Under the identification \eqref{eq:piv_character}, the ratio $j_{\DD}\circ (j|_{\DD})^{-1}$ determines a character $\chi_{\DD}\in\Hom(U(\DD),\BC^{\times})$.

\begin{lem}\label{lem:pivotal_structure_extension}
The following statements are equivalent:
\begin{enumerate}
    \item $j_{\DD}$ extends to a pivotal structure on $\CC$,
    \item $j_{\DD}$ and $j$ agree on $\DD\cap \CC_{\mathrm{ad}}$
    \item $\chi_{\DD}|_K=1$,
    \item $\chi_{\DD}$ factors through $H$. 
\end{enumerate}   
\end{lem}
\begin{proof}
$(1)\Rightarrow (2)$: Denote the extended pivotal structure by $j'$, $j,j'$ must be identical on $\CC_{\mathrm{ad}}=\CC_e$, which implies $(2)$.

$(2)\Rightarrow (3)$: Follows from the definition of $K$.

$(3)\Rightarrow (4)$: Follows from the short exact sequence.

$(4)\Rightarrow (1)$: Since $\CC$ is braided, $U(\CC)$ is abelian, character $\chi_{\DD}$ on $H\leq U(\CC)$ can be extended to a character $\chi\in\Hom(U(\CC),\BC^{\times})$. Therefore the modified pivotal structure $j'(X)=\chi(\deg(X)) j(X)$ extending $j_{\DD}$.
\end{proof}

Now combine Corollary \ref{cor:perturb_and_extension_center} and Lemma \ref{lem:pivotal_structure_extension}, we have the following corollary:
\begin{Cor}\label{cor:central_positivity_extension}
Let $\CC$ be a ribbon fusion category with corresponding spherical structure $j$, the following statements are equivalent:
\begin{enumerate}
    \item $d_j(X)>0$ for every nonzero object $X\in \CC'\cap\CC_{\mathrm{ad}}$,
    \item There exist an object $Z_v=X_{\mathrm{FP}}\otimes v\in(\CC')'_{\ZZ(\CC)}$, $v\in \Inv(\ZZ(\CC))$, such that the condition \eqref{eq:key_condition} holds.
\end{enumerate}
Note that condition (1) is independent of the pivotal structure on $\CC$.
\end{Cor}

\begin{rmk}\label{rmk:Galois_positive}
If $\CC$ satisfies the condition $d(X)>0$ for nonzero $X\in \CC'\cap\CC_{\mathrm{ad}}$, then every Galois conjugate
$\CC^\sigma$ satisfies the same condition, because these dimensions are integers, and taking M\"uger centers
and adjoint subcategories commutes with Galois conjugation.
\end{rmk}

\begin{Thm}\label{thm:main}
Let $\CC$ be a braided fusion category equipped with a pivotal structure $j$. Suppose that
\[
d_j^L(X)>0\qquad \text{for every nonzero }X\in \CC'\cap\CC_{\mathrm{ad}}.
\]
If $\CC$ has property F, then $\CC$ is weakly integral.  
\end{Thm}
\begin{proof}
We first prove the spherical case, the pivotal case will follow from sphericalization. For every Galois conjugate $\CC^{\sigma}$, from Remark \ref{rmk:Galois_positive} and Corollary \ref{cor:central_positivity_extension}, there always exists $Z_{\sigma}$ in $\operatorname{Im}(G)$, Property F of $\CC^{\sigma}$ induces the property F of $Z_{\sigma}$ through the dominant ribbon functor $G$. Therefore similar to the proof of Theorem \ref{thm:modular_case}, we have $\ZZ(\CC)^{\sigma}$ is pseudounitary for all $\sigma\in\Aut(\BC)$, which implies $\CC^{\sigma}$ is pseudounitary for all $\sigma\in\Aut(\BC)$, thus weakly integral. 

Now let $F:\widetilde{\CC} \longrightarrow\CC$ be the forgetful functor from the sphericalization, since $\CC$ is braided, its braiding lifts to $\widetilde{\CC}$, making $F$ a braided tensor functor. The pivotal structure $j$ lifts to $\widetilde{\CC}$ by setting $\widetilde j_{(X,\phi)}=j_X$, as $\phi^{**}j_X=j_{X^{**}}\phi=(j_X)^{**}\phi$. Therefore $\widetilde{j}$ and the canonical spherical structure $j^{\mathrm{sph}}$ agree on $\widetilde{\CC}_{\mathrm{ad}}$. The functor $F$ is essentially surjective, braided and sends simple objects to simple objects. Consequently,
\[
F\bigl(\widetilde{\CC}'\cap \widetilde{\CC}_{\mathrm{ad}}\bigr) \subseteq\CC'\cap\CC_{\mathrm{ad}},
\]
and for every simple object $Y\in \widetilde{\CC}'\cap \widetilde{\CC}_{\mathrm{ad}}$,
\[
d_{j^{\mathrm{sph}}}(Y)=d_j^L(F(Y))>0.
\]
Faithfulness and braidedness of $F$ also imply that $\widetilde{\CC}$ has property F whenever $\CC$ does. Now from the spherical case, $\widetilde{\CC}$ is weakly integral. Since
\[
\FPdim(\widetilde{\CC}) =2\FPdim(\CC),
\]
$\FPdim(\CC)$ is a rational algebraic integer, and hence an integer.
\end{proof}

\begin{Cor}
 If $\CC$ is a pseudounitary braided category and $X\in \Irr(\CC)$ has property F, then $\FPdim(X)^2\in\BZ$.
\end{Cor}
\begin{proof}
Property F and pseudounitarity both pass to the fusion subcategory generated by $X$, therefore, by Theorem \ref{thm:main}, the subcategory is weakly integral.
\end{proof}

\section{Positivity condition and the existence of pivotal
structures}\label{sec:equivalence}
In this section, we relate positivity of the pivotal structure induced on $\CC'\cap \CC_{\mathrm{ad}}$ to the existence of pivotal structures on braided fusion categories. We prove that this condition holds for every pivotal braided fusion category if and only if every braided fusion category over $\BC$ admits a pivotal structure. The idea of the proof is as follows:
\begin{itemize}
    \item Assuming positivity, the negative-dimensional transparent object $\varepsilon$ in the sphericalization $\widetilde{\CC}$ cannot lie in $\widetilde{\CC}_{\mathrm{ad}}$. A suitable grading character then modifies the canonical pivotal structure into one that descends to $\CC$.
    \item Conversely, we first reduce to the case when $\CC'$ is Tannakian (Lemma \ref{lem:superreduction}), then embed it into a central $\operatorname{Rep}(H)$ with centerless $H$ (Lemma \ref{lem:extension}). Then the unique pivotal structure on $\operatorname{Rep}(H)$ and agreement of pivotal structures on the adjoint subcategory forces the required positivity.
\end{itemize}

\begin{lem}\label{lem:centerless_cover}
Every finite group is a quotient of a finite centerless group $H$.
\end{lem}

\begin{proof}
Let $G$ be a finite group, and choose a nontrivial finite centerless group $K$, for example $K=S_3$. Consider the regular wreath product 
\[
H=K^G\rtimes G,
\] 
where $K^G$ is the group of functions $G\to K$ with pointwise multiplication, and $G$ acts on $K^G$ by left translation: $g\cdot f(x)=f(g^{-1}x)$. The projection $H\to G$ is surjective. It remains to show that $Z(H)=\{e\}$. Suppose $(f,g)\in Z(H)$. Since $(f,g)$ commutes with $(f',e)$ for every $f'\in K^G$, 
\begin{equation}\label{eq:wreath_proudct_center}
f(x)f'(g^{-1}x)=f'(x)f(x) \qquad (x\in G,\ f'\in K^{G}).   
\end{equation}
If $g\ne e$, choose $k\in K\setminus\{e\}$ and let $f'\in K^G$ satisfy $f'(e)=k$ and $f'(x)=e$ for $x\ne e$. Evaluating \eqref{eq:wreath_proudct_center} at $x=g$ yields $f(g)k=f(g)$, a contradiction. Thus $g=e$. Equation \eqref{eq:wreath_proudct_center} now implies that $f(x)$ commutes with every element of $K$, for every $x\in G$. Since $Z(K)=\{e\}$, we obtain $f(x)=e$ for all $x$. Hence $Z(H)=\{e\}$.
\end{proof}

\begin{lem}\label{lem:extension}
 Let $\CC$ be a braided fusion category with $\Rep(G)\subseteq \CC'$. For every surjective homomorphism $\pi: H\twoheadrightarrow G$ of finite groups, there exists a braided fusion category $\DD$ containing $\CC$ as a braided fusion subcategory, together with a Tannakian subcategory $\Rep(H)\subseteq \DD'$, such that 
 \[
\Rep(G)\subseteq\Rep(H)\subseteq\DD'.
 \]
\end{lem}
\begin{proof}
Since $\Rep(G)\subseteq \CC'$, de-equivariantization gives a braided fusion category $\CC_G$ equipped with an action of $G$ by braided tensor autoequivalences. Moreover, $\CC\simeq_{\mathrm{br}} (\CC_G)^G$ as braided fusion categories over $\Rep(G)$ \cite[Sec.~4]{DGNO}. Pulling the action back along $\pi$ gives an action of $H$ on $\CC_G$. Set $\DD=(\CC_G)^H$, which is again a braided fusion category. We show $\DD$ has the desired properties. 

Let $T_g$ denote the action of $G$, define
\[
I_\pi:(\CC_G)^G\longrightarrow(\CC_G)^H,\qquad\bigl(X,(u_g)_{g\in G}\bigr)\longmapsto\bigl(X,(u_{\pi(h)})_{h\in H}\bigr),
\]
and let $I_{\pi}$ act as identity on the underlying morphisms, the tensor products and the braidings inherited from $\CC_G$ makes $I_{\pi}$ a braided tensor functor. It is also fully faithful, a morphism $f:X\to Y$ in $\CC_G$ is a morphism in $\Hom_{(\CC_G)^H}(I_{\pi}(X,u),I_{\pi}(Y,v))$ precisely when 
\[
f\circ u_{\pi(h)}=v_{\pi(h)}\circ T_{\pi(h)}(f)\qquad \textit{for every $h\in H$}.
\]
Since $\pi$ is surjective, this is equivalent to
\[
f\circ u_g=v_g\circ T_g(f)\qquad \textit{for every $g\in G$}.
\]
Consequently, $f$ is a morphism in $\Hom_{\CC^G_G}((X,u),(Y,v))$. Thus $I_{\pi}$ is fully faithful.
\end{proof}

\begin{lem}\label{lem:superreduction}
Let $\CC$ be a pivotal braided fusion category. There exist a pivotal braided fusion category $\DD$, a Tannakian subcategory $\EE\subseteq\DD'$, and a pivotal tensor equivalence
\[
F:\CC'\xrightarrow{\sim}\EE
\]
such that
\[
F(\CC'\cap\CC_{\mathrm{ad}})=\EE\cap\DD_{\mathrm{ad}}.
\]
Here $\CC'$ and $\EE$ carry the pivotal structures inherited from $\CC$ and $\DD$, respectively.
\end{lem}
\begin{proof}
If $\CC'$ is Tannakian, take $\DD=\CC,\ \EE=\CC'$, and $F=\Id$. Otherwise Deligne's theorem gives a symmetric tensor equivalence $\CC'\simeq \Rep(G,z)$, where $z\in Z(G)$ and $z^2=1$, therefore under identification \eqref{eq:piv_character}, it induces a $\Zn2$ grading on $\Rep(G,z)$, the braiding for homogeneous objects $X, Y$ is given by 
\[
c_{X,Y}=(-1)^{\deg(X)\deg(Y)}\tau_{X,Y},
\]
where $\tau_{X,Y}$ is the flip.

Now consider the category $\BB=\PSU(2)_6$ with its standard unitary ribbon structure \cite[Sec.~3]{SMCat20}. It has simples $1,a,fa,f$ with $f$ a fermion in $\BB'$, the fusion rules $a^2=\one\oplus a\oplus fa$ and $a\otimes fa=f\oplus a\oplus fa$ put $f\in \BB_{\mathrm{ad}}$, therefore
\[
\quad c_{f,f}=-\Id,\quad d(f)=\FPdim(f)=1, \quad f\in \BB_\mathrm{ad}.
\]
Therefore, set
\[
\DD=:\CC\boxtimes \PSU(2)_6,
\]
with its product pivotal structure. 
Let $\EE$ be the fusion subcategory with simple objects
\[
X\boxtimes f^{\deg(X)}, \qquad X\in \Irr(\Rep(G,z))
\]
It is easy to see $\EE\simeq_{\mathrm{br}}\Rep(G)$ and $d^L_{\DD}(X\boxtimes f^{\deg(X)})=d^L_{\CC}(X)$, therefore $F: X\mapsto X\boxtimes f^{\deg(X)}$ extends to a pivotal tensor equivalence. The last statement follows from $\DD_{\mathrm{ad}}=\CC_{\mathrm{ad}}\boxtimes\BB_{\mathrm{ad}}$ and $f\in \BB_{\mathrm{ad}}$.
\end{proof}
\begin{Thm}
 The following two statements are equivalent
 \begin{enumerate}
     \item For every braided fusion category $\CC$, equipped with a pivotal structure $j$, 
     \[ 
     d^L_j(X)=\FPdim(X)\qquad \textit{for every nonzero}\ X\in\CC'\cap \CC_{\mathrm{ad}}.
     \]
     Equivalently, these left pivotal dimensions are positive.
     \item Every braided fusion category over $\BC$ admits a pivotal structure.
 \end{enumerate}  
\end{Thm}

\begin{proof}
$(1)\Rightarrow (2)$: Let $\CC$ be a braided fusion category and $F:\widetilde{\CC} \longrightarrow\CC$ be the forgetful functor from the sphericalization. The braiding of $\CC$ lifts to $\widetilde{\CC}$, thus the canonical spherical structure makes $\widetilde{\CC}$ a ribbon fusion category. Let $\varepsilon=(\one,-\operatorname{id}_{\one})$, using the canonical identification $\one^{**}\cong\one$. Then
\[
\varepsilon\in\widetilde{\CC}',\qquad\varepsilon^{\otimes2}\cong\one,\qquad d(\varepsilon)=-1,\qquad\FPdim(\varepsilon)=1.
\]
By (1), $\varepsilon\notin\widetilde{\CC}_{\mathrm{ad}}$.
Consequently, $t:=\deg(\varepsilon)$ is a nontrivial element of order two in $U(\widetilde{\CC})$. This group is finite abelian because $\widetilde{\CC}$ is
braided. Characters of a finite abelian group separate its elements, so there is a character
\[
\chi:U(\widetilde{\CC})\to\BC^\times\qquad\text{with}\qquad\chi(t)=-1.
\]
Modify the canonical pivotal structure by this character. For a simple lift $(X,\varphi_X)$, the resulting pivotal map has underlying morphism
\[
\widetilde j_{(X,\varphi_X)}
=\chi\bigl(\deg(X,\varphi_X)\bigr)\varphi_X.
\]
The two simple lifts of $X$ are $(X,\varphi_X)$ and $(X,-\varphi_X)\cong\varepsilon\otimes(X,\varphi_X)$.
Writing $g=\deg(X,\varphi_X)$, we obtain
\[
F\bigl(\widetilde j_{(X,-\varphi_X)}\bigr)
=\chi(tg)(-\varphi_X)
=\chi(g)\varphi_X
=F\bigl(\widetilde j_{(X,\varphi_X)}\bigr).
\]
Define $j_X: X\to X^{**}$ for each simple $X$ by this common underlying morphism and extend by semisimplicity. The resulting natural isomorphism is monoidal: choose lifts of two simple objects, decompose their tensor product in $\widetilde{\CC}$, and apply monoidality of $\widetilde j$ and the independence just proved. Thus $j$ is a pivotal structure on $\CC$.

$(2)\Rightarrow (1)$: Let $\CC$ be a pivotal braided fusion category. By Lemma \ref{lem:superreduction}, choose a pivotal braided fusion category $\DD_0$, a Tannakian subcategory $\EE\subseteq\DD_0'$, and a pivotal tensor equivalence $F:\CC'\to\EE$ satisfying
\[
F(\CC'\cap\CC_{\mathrm{ad}})=\EE\cap(\DD_0)_{\mathrm{ad}}.
\]
Identify $\EE$ with $\Rep(G)$ as a braided fusion category. By Lemma \ref{lem:centerless_cover}, choose a surjection $H\twoheadrightarrow G$ with $Z(H)=\{e\}$. Lemma \ref{lem:extension} gives a braided fusion category $\DD$ containing $\DD_0$ such that
\[
\EE\simeq\Rep(G)\subseteq\Rep(H)\subseteq\DD'.
\]
By (2), $\DD$ admits a pivotal structure $j$. Let $j_0$ denote the given pivotal structure on $\DD_0$. We have
\[
j|_{(\DD_0)_{\mathrm{ad}}}=j_0|_{(\DD_0)_{\mathrm{ad}}}.
\]
On the other hand,
\[
\Aut_{\otimes}\bigl(\Id_{\Rep(H)}\bigr)\cong Z(H)=\{e\}.
\]
Hence $\Rep(H)$ has a unique pivotal structure, its
canonical positive one. The restriction of $j$ to $\Rep(H)$ must be this structure. For every simple $Y\in\EE\cap(\DD_0)_{\mathrm{ad}}$, it follows that
\[
d^L_{j_0}(Y)=d^L_j(Y)=\FPdim(Y).
\]
Since $F$ is pivotal and preserves Frobenius--Perron dimensions, the required identity follows for every simple $X\in\CC'\cap\CC_{\mathrm{ad}}$.

\end{proof}

\section*{Data Availability}
Data sharing is not applicable to this article as no datasets were generated or analyzed during the current study.

\section*{Conflict of interest}
The authors have no conflicts of interest to declare that are relevant to the content of this article.

\bibliographystyle{abbrv}
\bibliography{Reference}

@article{Rowell2010,
  author  = {Rowell, Eric C.},
  title   = {Braid representations from quantum groups of exceptional {Lie} type},
  journal = {Revista de la Uni{\'o}n Matem{\'a}tica Argentina},
  volume  = {51},
  number  = {1},
  pages   = {165--175},
  year    = {2010},
  url     = {https://inmabb.criba.edu.ar/revuma/pdf/v51n1/v51n1a08.pdf}
}

@article{RW2012,
  author  = {Rowell, Eric C. and Wang, Zhenghan},
  title   = {Localization of unitary braid group representations},
  journal = {Communications in Mathematical Physics},
  volume  = {311},
  number  = {3},
  pages   = {595--615},
  year    = {2012},
  doi     = {10.1007/s00220-011-1386-7},
  url     = {https://web.math.ucsb.edu/~zhenghwa/data/research/pub/Localization-12.pdf}
}

@article{GHR2013,
  author  = {Galindo, C{\'e}sar and Hong, Seung-Moon and Rowell, Eric C.},
  title   = {Generalized and quasi-localizations of braid group representations},
  journal = {International Mathematics Research Notices},
  volume  = {2013},
  number  = {3},
  pages   = {693--731},
  year    = {2013},
  url     = {https://arxiv.org/abs/1105.5048}
}

@Book{BK01,
  author     = {Bakalov, Bojko and Kirillov, Jr., Alexander},
  title      = {Lectures on tensor categories and modular functors},
  isbn       = {0-8218-2686-7},
  pages      = {x+221},
  publisher  = {American Mathematical Society, Providence, RI},
  series     = {University Lecture Series},
  volume     = {21},
  mrclass    = {18D10 (17B37 57M27 57R56 81R50 81T40 81T45)},
  mrnumber   = {1797619},
  mrreviewer = {J. Stasheff},
  year       = {2001},
}

@article {SMCat20,
    AUTHOR = {Bruillard, Paul and Galindo, C\'{e}sar and Ng, Siu-Hung and
              Plavnik, Julia Y. and Rowell, Eric C. and Wang, Zhenghan},
     TITLE = {Classification of super-modular categories by rank},
   JOURNAL = {Algebr. Represent. Theory},
  FJOURNAL = {Algebras and Representation Theory},
    VOLUME = {23},
      YEAR = {2020},
    NUMBER = {3},
     PAGES = {795--809},
      ISSN = {1386-923X,1572-9079},
   MRCLASS = {18M20 (16D90 81U20)},
  MRNUMBER = {4109138},
MRREVIEWER = {Costel\ Gabriel\ Bontea},
       DOI = {10.1007/s10468-019-09873-9},
       URL = {https://doi.org/10.1007/s10468-019-09873-9},
}

@Article{Del02,
  author     = {Deligne, P.},
  title      = {Cat\'{e}gories tensorielles},
  doi        = {10.17323/1609-4514-2002-2-2-227-248},
  issn       = {1609-3321},
  number     = {2},
  pages      = {227--248},
  url        = {https://doi.org/10.17323/1609-4514-2002-2-2-227-248},
  volume     = {2},
  fjournal   = {Moscow Mathematical Journal},
  journal    = {Mosc. Math. J.},
  mrclass    = {18D10 (18E30 20C99)},
  mrnumber   = {1944506},
  mrreviewer = {Viktor Ostrik},
  year       = {2002},
}

@Book{EGNO,
  author     = {Etingof, Pavel and Gelaki, Shlomo and Nikshych, Dmitri and Ostrik, Victor},
  title      = {Tensor categories},
  doi        = {10.1090/surv/205},
  isbn       = {978-1-4704-2024-6},
  pages      = {xvi+343},
  publisher  = {American Mathematical Society, Providence, RI},
  series     = {Mathematical Surveys and Monographs},
  url        = {http://dx.doi.org/10.1090/surv/205},
  volume     = {205},
  mrclass    = {18D10 (16T05)},
  mrnumber   = {3242743},
  mrreviewer = {Julien Bichon},
  year       = {2015},
}

@Article{ENO05,
  author     = {Etingof, Pavel and Nikshych, Dmitri and Ostrik, Viktor},
  title      = {On fusion categories},
  doi        = {10.4007/annals.2005.162.581},
  issn       = {0003-486X},
  number     = {2},
  pages      = {581--642},
  url        = {http://dx.doi.org/10.4007/annals.2005.162.581},
  volume     = {162},
  fjournal   = {Annals of Mathematics. Second Series},
  journal    = {Ann. of Math. (2)},
  mrclass    = {16W30 (18D10)},
  mrnumber   = {2183279},
  mrreviewer = {Yorck Sommerh\"auser},
  year       = {2005},
}

@Article{GN08,
  author     = {Gelaki, Shlomo and Nikshych, Dmitri},
  title      = {Nilpotent fusion categories},
  doi        = {10.1016/j.aim.2007.08.001},
  issn       = {0001-8708},
  number     = {3},
  pages      = {1053--1071},
  url        = {https://doi-org.libezp.lib.lsu.edu/10.1016/j.aim.2007.08.001},
  volume     = {217},
  fjournal   = {Advances in Mathematics},
  journal    = {Adv. Math.},
  mrclass    = {18D99 (20D15)},
  mrnumber   = {2383894},
  mrreviewer = {I. B. S. Passi},
  year       = {2008},
}

@article {NR11,
    AUTHOR = {Naidu, Deepak and Rowell, Eric C.},
     TITLE = {A finiteness property for braided fusion categories},
   JOURNAL = {Algebr. Represent. Theory},
  FJOURNAL = {Algebras and Representation Theory},
    VOLUME = {14},
      YEAR = {2011},
    NUMBER = {5},
     PAGES = {837--855},
      ISSN = {1386-923X},
   MRCLASS = {18D10 (20F36)},
  MRNUMBER = {2832261},
MRREVIEWER = {Alexei Davydov},
       DOI = {10.1007/s10468-010-9219-5},
       URL = {http://dx.doi.org.proxy.lib.ohio-state.edu/10.1007/s10468-010-9219-5},
}

@Article{RSW09,
  author     = {Rowell, Eric and Stong, Richard and Wang, Zhenghan},
  title      = {On classification of modular tensor categories},
  doi        = {10.1007/s00220-009-0908-z},
  issn       = {0010-3616},
  number     = {2},
  pages      = {343--389},
  url        = {https://doi-org.libezp.lib.lsu.edu/10.1007/s00220-009-0908-z},
  volume     = {292},
  fjournal   = {Communications in Mathematical Physics},
  journal    = {Comm. Math. Phys.},
  mrclass    = {18D10 (81P68 81R05)},
  mrnumber   = {2544735},
  mrreviewer = {Ram\'{o}n Gonz\'{a}lez Rodr\'{\i}guez},
  year       = {2009},
}

@article{RW18,
	AUTHOR = {Rowell, Eric C. and Wang, Zhenghan},
	TITLE = {Mathematics of topological quantum computing},
	JOURNAL = {Bull. Amer. Math. Soc. (N.S.)},
	FJOURNAL = {American Mathematical Society. Bulletin. New Series},
	VOLUME = {55},
	YEAR = {2018},
	NUMBER = {2},
	PAGES = {183--238},
	ISSN = {0273-0979},
	MRCLASS = {81P68 (03D10 18D10 57R56 68Q12 81T45)},
	MRNUMBER = {3777017},
	DOI = {10.1090/bull/1605},
	URL = {https://doi-org.libezp.lib.lsu.edu/10.1090/bull/1605},
}

@article {DGNO,
    AUTHOR = {Drinfeld, Vladimir and Gelaki, Shlomo and Nikshych, Dmitri and
              Ostrik, Victor},
     TITLE = {On braided fusion categories. {I}},
   JOURNAL = {Selecta Math. (N.S.)},
  FJOURNAL = {Selecta Mathematica. New Series},
    VOLUME = {16},
      YEAR = {2010},
    NUMBER = {1},
     PAGES = {1--119},
      ISSN = {1022-1824},
   MRCLASS = {18D10 (16T05)},
  MRNUMBER = {2609644},
MRREVIEWER = {R. Gylys},
       DOI = {10.1007/s00029-010-0017-z},
       URL = {https://doi-org.srv-proxy1.library.tamu.edu/10.1007/s00029-010-0017-z},
}

@Book{Tur94,
  author     = {Turaev, V. G.},
  title      = {Quantum invariants of knots and 3-manifolds},
  isbn       = {3-11-013704-6},
  pages      = {x+588},
  publisher  = {Walter de Gruyter \& Co., Berlin},
  series     = {de Gruyter Studies in Mathematics},
  volume     = {18},
  mrclass    = {57M25 (16W30 17B37 57N10 81R50)},
  mrnumber   = {1292673},
  mrreviewer = {Louis H. Kauffman},
  year       = {1994},
}

@article {GRR20,
    AUTHOR = {Gustafson, Paul and Rowell, Eric C. and Ruan, Yuze},
     TITLE = {Metaplectic categories, gauging and property {$F$}},
   JOURNAL = {Tohoku Math. J. (2)},
  FJOURNAL = {The Tohoku Mathematical Journal. Second Series},
    VOLUME = {72},
      YEAR = {2020},
    NUMBER = {3},
     PAGES = {411--424},
      ISSN = {0040-8735},
   MRCLASS = {18M15 (20F36)},
  MRNUMBER = {4154826},
MRREVIEWER = {Costel Gabriel Bontea},
       DOI = {10.2748/tmj/1601085623},
       URL = {https://doi.org/10.2748/tmj/1601085623},
}

@article {GN21,
    AUTHOR = {Green, Jason and Nikshych, Dmitri},
     TITLE = {On the braid group representations coming from weakly
              group-theoretical fusion categories},
   JOURNAL = {J. Algebra Appl.},
  FJOURNAL = {Journal of Algebra and its Applications},
    VOLUME = {20},
      YEAR = {2021},
    NUMBER = {1},
     PAGES = {Paper No. 2150210, 20},
      ISSN = {0219-4988},
   MRCLASS = {18M05 (20F36)},
  MRNUMBER = {4209972},
MRREVIEWER = {Andr\'{e} G. Henriques},
       DOI = {10.1142/S0219498821502108},
       URL = {https://doi.org/10.1142/S0219498821502108},
}

@article {RW17,
    AUTHOR = {Rowell, Eric C. and Wenzl, Hans},
     TITLE = {{${\rm SO}(N)_2$} braid group representations are {G}aussian},
   JOURNAL = {Quantum Topol.},
  FJOURNAL = {Quantum Topology},
    VOLUME = {8},
      YEAR = {2017},
    NUMBER = {1},
     PAGES = {1--33},
      ISSN = {1663-487X},
   MRCLASS = {20F36 (17B37 18D10 20G43)},
  MRNUMBER = {3630280},
MRREVIEWER = {Tim Raeymaekers},
       DOI = {10.4171/QT/85},
       URL = {https://doi.org/10.4171/QT/85},
}

@article {ERW08,
    AUTHOR = {Etingof, Pavel and Rowell, Eric and Witherspoon, Sarah},
     TITLE = {Braid group representations from twisted quantum doubles of
              finite groups},
   JOURNAL = {Pacific J. Math.},
  FJOURNAL = {Pacific Journal of Mathematics},
    VOLUME = {234},
      YEAR = {2008},
    NUMBER = {1},
     PAGES = {33--41},
      ISSN = {0030-8730,1945-5844},
   MRCLASS = {20F36 (16W30)},
  MRNUMBER = {2375313},
MRREVIEWER = {Stephen\ P.\ Humphries},
       DOI = {10.2140/pjm.2008.234.33},
       URL = {https://doi.org/10.2140/pjm.2008.234.33},
}

@article{Etingof2002Vafa,
  author        = {Etingof, Pavel},
  title         = {On {Vafa}'s theorem for tensor categories},
  journal       = {Mathematical Research Letters},
  volume        = {9},
  number        = {5--6},
  pages         = {651--657},
  year          = {2002},
  doi           = {10.4310/MRL.2002.v9.n5.a8},
  eprint        = {math/0207007},
  archivePrefix = {arXiv},
  primaryClass  = {math.QA}
}

@article{Rowell2005,
  author  = {Rowell, Eric C.},
  title   = {On a family of non-unitarizable ribbon categories},
  journal = {Mathematische Zeitschrift},
  volume  = {250},
  number  = {4},
  pages   = {745--774},
  year    = {2005},
  doi     = {10.1007/s00209-005-0773-1},
  url     = {https://people.tamu.edu/~rowell/RibbonB.pdf}
}

@misc{LiRuan2026,
  author        = {Enyu Li and Yuze Ruan},
  title         = {Generalized property {F} for weakly group-theoretical modular fusion categories},
  year          = {2026},
  note          = {arXiv:2610.05177},
  eprint        = {2610.05177},
  archivePrefix = {arXiv},
  primaryClass  = {math.QA},
  doi           = {10.48550/arXiv.2610.05177},
  url           = {https://arxiv.org/abs/2610.05177}
}

@article {Ostrik15,
    AUTHOR = {Ostrik, Victor},
     TITLE = {Pivotal fusion categories of rank 3},
   JOURNAL = {Mosc. Math. J.},
  FJOURNAL = {Moscow Mathematical Journal},
    VOLUME = {15},
      YEAR = {2015},
    NUMBER = {2},
     PAGES = {373--396, 405},
      ISSN = {1609-3321,1609-4514},
   MRCLASS = {18D10},
  MRNUMBER = {3427429},
MRREVIEWER = {Ehud\ Meir},
       DOI = {10.17323/1609-4514-2015-15-2-373-396},
       URL = {https://doi.org/10.17323/1609-4514-2015-15-2-373-396},
}
\end{document}